\documentclass[11pt]{article}

\usepackage[sort]{cite}       
\usepackage{amssymb}          

\usepackage{latexsym}         
\usepackage{amsmath}          
\usepackage[latin1]{inputenc} 
\usepackage[T1]{fontenc}
\usepackage{eucal}            
\usepackage{bbm}              
\usepackage{theorem}          
\usepackage{stmaryrd}         

\usepackage[english]{babel}
\usepackage{color}

\title{Formal Deformations of Dirac Structures}

\author{\textbf{Frank Keller}\thanks{E-mail: Frank.Keller@uni.lu}
  \\[0.1cm]
  University of Luxembourg\\
  Mathematics Laboratory\\
  162A, avenue de la Fa{\"{\i}}encerie\\
  L-1511 Luxembourg City\\
  Luxembourg\\[0.5cm]
  \addtocounter{footnote}{5}
  \textbf{Stefan Waldmann}\thanks{E-mail: Stefan.Waldmann@physik.uni-freiburg.de}
  \\[0.1cm]
  Fakult{\"a}t f{\"u}r Mathematik und Physik\\
  Albert-Ludwigs-Universit{\"a}t Freiburg\\
  Physikalisches Institut\\
  Hermann Herder Stra{\ss}e 3\\
  D 79104 Freiburg\\
  Germany}

\date{June 2006\\[0.5cm] FR-THEP 2006/09}

\renewcommand{\mathbb}[1]{\mathbbm{#1}} 

\newcounter{comment}

\newcommand{\id}         {\operatorname{\mathsf{id}}}

\newcommand{\Lie}        {\operatorname{\mathcal{L}}}    
  
\newcommand{\graph}      {\operatorname{\mathrm{graph}}} 
\newcommand{\Pol}        {\operatorname{\mathrm{Pol}}}

\newcommand{\End}        {\operatorname{\mathsf{End}}}

\newcommand{\varsharp}   {{\scriptscriptstyle \#}}

\newcommand{\ra}{\longrightarrow}
\newcommand{\R}{\mathbb{R}}
\newcommand{\F}{\mathcal F}
\newcommand{\B}{\mathcal B}

\newcommand{\scal}[1]{\langle{#1}\rangle}
\newcommand{\cour}[1]{[#1]_{{}_\mathcal{C}}}
\newcommand{\bra}[2]{{[#2]}_{\scriptscriptstyle \! #1}}
\newcommand{\roth}[1]{\{#1\}_{{}_{\!\mathcal{R}}}}

\newcommand{\der}[2]{[#2]_{{}_{\!#1}}}
\newcommand{\derT}[1]{[#1]_{{}_{\Theta}}}
\newcommand{\bili}[1]{(#1)}

\newcommand{\anc}[1]{\rho_{\scriptscriptstyle{#1}}}

\newcommand{\D}{\mathcal{D}}
\newcommand{\DT}{\D_{\scriptscriptstyle\Theta}}

\newcommand{\dcd}{\,\cdot\;\!,\:\!\cdot\,}

\newcommand{\op}{\operatorname}
\renewcommand{\(}{\left(}
\renewcommand{\)}{\right)}

\newcommand{\ta}{\tau^\varsharp a}
\newcommand{\ts}{\tau^\varsharp s}
\newcommand{\tu}{\tau^\varsharp u}

\newcommand{\coordup}[2]{#1^1,\ldots,#1^{#2}}
\newcommand{\coorddown}[2]{#1_1,\ldots,#1_{#2}}

\newcommand{\stCA}[1][M]{T #1 \oplus T^\ast\! #1}
\newcommand{\dif}{\operatorname{d}}
\newcommand{\dL}{\dif_L}
\newcommand{\partdif}[1]{\frac{\partial}{\partial #1}}

\newcommand{\sect}{\Gamma^\infty}
\newcommand{\Wedge}{\textstyle{\bigwedge}}
\newcommand{\pol}{\mathcal{P}}
\newcommand{\degL}{\deg_L}
\newcommand{\degLs}{\deg_{L^\ast}}
\newcommand{\grassalg}[1]{\sect(\Wedge^\bullet #1)}
\newcommand{\spalg}{\sect\big(\Wedge^\bullet\tau^\varsharp(L\oplus L^\ast)\big)}
\newcommand{\ldot}{\,\cdot\,}
\newtheorem{lemma}{Lemma}[section]
\newtheorem{proposition} [lemma] {Proposition}

\newtheorem{theorem} [lemma] {Theorem}
\newtheorem{corollary} [lemma] {Corollary}
\newtheorem{definition}[lemma] {Definition}

\theorembodyfont{\rmfamily}
\newtheorem{example}[lemma]{Example}
\newtheorem{remark}[lemma]{Remark}

\newenvironment{proof}[1][{}]{
  \par\noindent
  \textsc{Proof{#1}:}
}
{
 \hspace*{\fill}  $\blacksquare$\medskip\par 
}

\numberwithin{equation}{section}

\begin{document}

\maketitle

\begin{abstract}
    In this paper we set-up a general framework for a formal
    deformation theory of Dirac structures. We give a parameterization
    of formal deformations in terms of two-forms obeying a cubic
    equation. The notion of equivalence is discussed in detail. We
    show that the obstruction for the construction of deformations
    order by order lies in the third Lie algebroid cohomology of the
    Dirac structure. However, the classification of inequivalent first
    order deformations is not given by the second Lie algebroid
    cohomology but turns out to be more complicated.
\end{abstract}

\newpage

%
%

\tableofcontents

%
%

\section{Introduction}
\label{sec:Intro}

A Dirac structure is a maximally isotropic subbundle of a Courant
algebroid whose sections in addition are closed under the Courant
bracket. A Courant algebroid is a vector bundle with a not necessarily
positive definite fiber metric over a base manifold which is equipped
with a bundle map into the tangent bundle (the anchor) and a bracket
on its sections, the Courant bracket, subject to certain compatibility
conditions. The fundamental example of a Courant algebroid is $E = TM
\oplus T^*M$ with the natural pairing as fiber metric, the identity on
the first component as anchor and the bracket
\[
\cour{(X, \alpha), (Y, \beta)} 
= \left([X, Y], \Lie_X \beta - i_Y \dif \alpha\right).
\]
Then both $TM$ and $T^*M$ are Dirac structures in this Courant
algebroid.

Dirac structures were introduced by Courant \cite{courant:1990a} to
generalize on one hand symplectic and Poisson structures, on the other
hand they provide powerful tools to describe dynamics subject to
constraints. Moreover, they can also be used to encode various
`oid'-structures, in particular Lie bialgebroids
\cite{liu.weinstein.xu:1997a, mackenzie.xu:1994a, mackenzie:2005a}.
For the general notions of Dirac structures we refer to
\cite{courant:1990a, liu.weinstein.xu:1997a}.

As Dirac structures combine symplectic and Poisson structures and can
be used in constraint dynamics, it is natural to ask what a physically
reasonable quantization of a Dirac structure should be. In particular,
this could shine some new light on the quantization of constraint
dynamics and phase space reduction.  In deformation quantization
\cite{bayen.et.al:1978a}, see e.g.~\cite{weinstein:1994a,
  dito.sternheimer:2002a, gutt:2000a} for recent reviews, one knows
that the equivalence classes of formal quantizations of Poisson
structures are in one-to-one correspondence with equivalence classes
of formal deformations of the given Poisson structure into formal
Poisson structures modulo formal diffeomorphisms. This is one of the
main corollaries of Kontsevich's formality theorem
\cite{kontsevich:2003a}.

Motivated by this result, we investigate the deformation theory of
Dirac structures in order to determine their \emph{classical
  deformations} into formal Dirac structures up to formal
diffeomorphisms. We hope that this gives eventually some hints on how
to formulate a definition of deformation quantization of Dirac
structures such that a type of formality might hold true also in this
context. First steps in this direction have been taken in
\cite{severa:2005a:pre} by \v{S}evera who proposed a deformation
quantization of formal deformations of regular Dirac structures.  Note
however, that the classical deformations are also of interest if one
wants to describe stability/rigidity of Dirac structures, not
necessarily aiming at quantization. Thus our first aim of this paper
is to set-up a reasonable definition of a formal Dirac structure and
investigate basic properties of the corresponding classical
deformation theory.

As a Dirac structure $L$ gives in particular the structure of a Lie
algebroid it is natural to compare the formal deformation theory of
Dirac structures with the formal deformation theory of $L$ as a Lie
algebroid in the sense of \cite{crainic.moerdijk:2004a:pre}: It turns out
that any deformation of a Dirac structure induces a Lie
algebroid deformation, but not necessarily vice-versa. Moreover, in
\cite{crainic.moerdijk:2004a:pre} it was shown that the Lie algebroid
structure of $TM$ is rigid with respect to formal deformations while
it is easy to see that this is not the case for Dirac structure
deformations, here any non-trivial pre-symplectic form provides a
non-trivial deformation.

The main results of this paper is on one hand, that the obstruction
space for formal order-by-order deformations of a Dirac structure $L$
is given by the third Lie algebroid cohomology of $L$. On the other
hand, and this is the more surprising result, the reasonable notion of
equivalence up to formal diffeomorphisms does \emph{not} yield a
classification of inequivalent first order deformations in the second
Lie algebroid cohomology, as one might first think: the actual
classification is more involved and seems to be beyond a simple
cohomological formulation. This depends of course on our definition of
equivalence which we based on formal diffeomorphisms. Most of our
results emerged from the Diplomarbeit \cite{keller:2004a}.

As main technique it turned out that a description of formal Dirac
structures in terms of graphs of formal two-forms requires some
reasonable calculus. We found the derived bracket formalism
\cite{kosmann-schwarzbach:1996a,kosmann-schwarzbach:2004b}, already
introduced by Roytenberg in a super-geometric way
\cite{roytenberg:1999a}, most useful. However, we realized the derived
bracket formalism not in terms of super-geometry but used more
conventional geometric objects: the main ingredient is the
Rothstein-Poisson bracket \cite{rothstein:1991a}. We believe that this
approach has its own interest, in particular when it will come to
quantization as we can rely on Bordemann's results for the deformation
quantization of the Rothstein-Poisson bracket
\cite{bordemann:1996a:pre, bordemann:2000a}. Nevertheless, our
approach is completely equivalent to the one of Roytenberg.

The paper is organized as follows: In
Section~\ref{sec:GeneralDiracStuff} we recall some basic definitions
and results on Courant algebroids, their automorphisms and their Dirac
structures. Section~\ref{sec:DerivedBrackets} introduces the derived
bracket point of view in order to handle the quite complicated
algebraic identities of the Courant bracket in a more efficient way.
We recall the Rothstein-Poisson bracket and use it to formulate Dirac
structures in this context. In
Section~\ref{sec:SmoothFormalDeformations} we first formulate a smooth
deformation of a Dirac structure and discuss the problem of
equivalence up to diffeomorphisms for the case of a general Courant
algebroid. Taking this as motivation we pass to formal deformations by
Taylor expansion in the deformation parameter as usual. The
fundamental equation, a sort of Maurer-Cartan equation which controls
the deformation, has already been discussed in some different context
in \cite[Eq.~(4.3)]{roytenberg:2002a}. We show that the order-by-order
construction of a formal deformation yields obstructions in the third
Lie algebroid cohomology of the undeformed Dirac structure. Finally,
we discuss the notion of equivalence up to formal diffeomorphisms in
detail and point out that the second Lie algebroid cohomology is not
necessarily the space of inequivalent first order deformations.
Finally, Appendix~\ref{sec:Rothstein} gives an overview on the
Rothstein-Poisson bracket and recalls some of its basic properties.

\smallskip

\noindent
\textbf{Conventions:} Throughout the paper we use Einstein's summation
convention, i.e. summation over repeated coordinate indices is
automatic.

\medskip

\noindent
\textbf{Acknowledgments:} We would like to thank Lorenz Schwachhöfer
for valuable discussions on the notion of equivalent deformations and
Marco Gualtieri for a stimulating remark on the quantization aspect.
Moreover, we would like to thank Pavol \v{S}evera, Jim Stasheff and
Alan Weinstein for valuable remarks and suggestions on the first
version.

%
%

\section{General Remarks on Dirac Structures in Courant Algebroids}
\label{sec:GeneralDiracStuff}

In this section we recall some basic notions of Courant algebroids and
Dirac structures in order to set up our notation. Most of the material
is standard, see e.g.~\cite{courant:1990a, liu.weinstein.xu:1997a,
  roytenberg:1999a}.

%
%

\subsection{Courant Algebroids}
\label{subsec:CourantAlgebroids}

\begin{definition}
    \label{def.courant-algebroid}
    A Courant algebroid is a vector bundle $E \ra M$ together with
    nondegenerate symmetric bilinear form $h$, a bracket $\cour{\dcd}:
    \sect(E) \times \sect(E) \ra \sect(E)$ on the sections of the
    bundle and a vector bundle homomorphism $\rho: E \ra TM$, called
    anchor, such that for all $e_1,e_2,e_3 \in \sect(E)$ and $f\in
    C^\infty(M)$ the following conditions hold:
    \begin{enumerate}
    \item Jacobi identity, i.e.  $\cour{e_1, \cour{e_2,
            e_3}}=\cour{\cour{e_1, e_2}, e_3} + \cour{e_2, \cour{e_1,
            e_3}},$
    \item $\cour{e_1, e_2} + \cour{e_2, e_1}= \D\, h(e_1, e_2)$, where
        $\D: C^\infty(M) \ra \sect(E)$ is defined by
        $$h(\D f, e) = \rho(e) f,$$
    \item $\rho(e_1)h(e_2, e_3)= h(\cour{e_1, e_2}, e_3) + h(e_2,
        \cour{e_1, e_3})$.
\end{enumerate}
\end{definition}
An easy computation shows that the Courant bracket $\cour{\cdot,
  \cdot}$ satisfies the Leibniz rule
\begin{equation}
    \label{eq:LeibnizRule}
    \cour{e_1, f e_2} = f \cour{e_1, e_2} + (\rho(e_1) f) e_2
\end{equation}
and the anchor turns out to satisfy
\begin{equation}
    \label{eq:AnkerLieMorph}
    \rho(\cour{e_1, e_2}) = 
    [\rho(e_1),\rho(e_2)]
\end{equation}
for all $e_1, e_2, e_3 \in \Gamma^\infty(E)$ and $f \in C^\infty(M)$,
see e.g. \cite{uchino:2002a, kosmann-schwarzbach:2005a,
  kosmann-schwarzbach.magri:1990a}.
\begin{remark}
    Equivalent to this definition is the one given in
    \cite{carinena.grabowski.marmo:2004a}. One can also consider the
    object obtained by skew-symmetrization of the Courant bracket,
    which is sometimes referred to as a Courant algebroid. Both
    definitions are equivalent, see \cite{roytenberg:1999a} for a
    detailed discussion.
\end{remark}

The above definition for a Courant algebroid is the generalization of
an object studied by Courant in \cite{courant:1990a}, which we will
refer as the standard Courant algebroid:
\begin{example}[Standard Courant algebroid \cite{courant:1990a}]
    Consider the vector bundle $E= \stCA$ over a manifold $M$. The
    canonical symmetric bilinear form on $E$ given by
    \begin{equation}
        \label{eq:NaturalPairing}
        \scal{(X,\alpha),(Y,\beta)} = \alpha(Y) + \beta(X),
    \end{equation}
    where $X, Y \in \mathfrak{X}(M)$ and $\alpha,\beta \in
    \Omega^1(M)$, together with the bracket
    \begin{equation}
        \label{eq:StandardCourantBracket}
        \cour{(X,\alpha),(Y,\beta)} = 
        ([X,Y],\Lie_X \beta - i_Y \dif \alpha)
    \end{equation}
    and the anchor $\rho$ defined by $\rho(X,\alpha) = X$ endows $E$
    with the structure of a Courant algebroid.
\end{example}
\begin{remark}
    According to our definition of a Courant algebroid we use here
    also the non skew-symmetric version where originally in
    \cite{courant:1990a} the skew-symmetric bracket was used.
\end{remark}
Other examples for Courant algebroids are given by the double of Lie
bialgebroids \cite{liu.weinstein.xu:1997a}, or more general by the
doubles of Lie quasi-bialgebroids or proto bialgebroids, see
\cite{kosmann-schwarzbach:2005a}.  We will come back to these examples
later.

%
%

\subsection{Automorphisms of Courant Algebroids}
\label{subsec:Automorphisms}

Crucial for our investigations of deformations of Dirac structures
will be an appropriate notion of isomorphism. To this end we need the
automorphisms of the Courant algebroid. If $E \longrightarrow M$ is a
Courant algebroid, then a vector bundle automorphisms $\Phi: E
\longrightarrow E$ over a diffeomorphism $\phi: M \longrightarrow M$
is called an automorphism of the Courant algebroid, if the following
two conditions are fulfilled: First, $\Phi$ is an isometry of the
bilinear form $h$, i.e. for all $e_1, e_2 \in \sect(E)$
\begin{equation}
    \label{eq:Isometry}
    h(\Phi^\ast e_1, \Phi^\ast e_2) = \phi^*(h(e_1, e_2)).
\end{equation}
Second, $\Phi$ is natural with respect to the Courant bracket, i.e.
for all $e_1,e_2 \in \sect(E)$
\begin{equation}
    \label{eq:PhiPreservesBracket}
    \cour{\Phi^\ast e_1, \Phi^\ast e_2} = \Phi^\ast \cour{e_1, e_2}.
\end{equation}
The following lemma shows that the compatibility with the anchor is
already fixed by these two conditions:
\begin{lemma}
    \label{lemma.couralgauto}
    If $\Phi:E \ra E$ is a Courant algebroid automorphism then  the
    anchor  $\rho$ satisfies
    \begin{equation}
        \label{eq:AutomorphismAnchor}
        \rho \circ \Phi = T\phi \circ \rho.
    \end{equation}
\end{lemma}
\begin{proof} This is used implicitly in
    \cite[Prop.~3.24]{gualtieri:2003a}: Using \eqref{eq:LeibnizRule}
    and then \eqref{eq:PhiPreservesBracket} gives $\cour{\Phi^\ast
      e_1, \Phi^\ast (f e_2)} = \Phi^\ast(f \cour{e_1,e_2}) +
    \rho(\Phi^\ast e_1)(\phi^\ast f) \Phi^\ast e_2$.  The other way
    round gives $\cour{\Phi^\ast e_1, \Phi^\ast (f e_2)} = \Phi^\ast(f
    [e_1,e_2]) + \phi^\ast(\rho(e_1) f)\Phi^\ast e_2$ whence we obtain
    \[
    \rho(\Phi^\ast e_1)(\phi^\ast f) 
    = \phi^\ast(\rho(e_1)f) 
    =  \phi^\ast(\rho(e_1))(\phi^\ast f)
    \]
    for all $e_1, e_2 \in \Gamma^\infty(E)$ and $f \in C^\infty(M)$,
    which implies \eqref{eq:AutomorphismAnchor}.
\end{proof}

In case of the standard Courant algebroid one can determine the group
of automorphisms completely. We recall the following definition
\cite{severa.weinstein:2001a}:
\begin{definition}[Gauge Transformations]
    \label{definition:GaugeTrafo}
    Let $E = \stCA$ be the standard Courant algebroid and $B \in
    \Omega^2(M)$ a two-form. A gauge transformation is a map $\tau_B:
    \stCA \ra \stCA$ given by $\tau_B(X,\alpha) = (X,\alpha + i_X B)$.
\end{definition}
\begin{lemma}[\v{S}evera, Weinstein \cite{severa.weinstein:2001a}]
    \label{lemma.gaugetrans}
    A gauge transformation $\tau_B$ is an automorphism of the standard
    Courant algebroid structure on $TM\oplus T^\ast M$ if and only if
    $B$ is closed.
\end{lemma}

Let $\phi$ be a diffeomorphism of $M$. Then we denote the canonical
lift of $\phi$ to $\stCA$ by $\F \phi = (T\phi,T_\ast \phi)$, where
$T_\ast \phi: T^\ast M\ra T^\ast M$ is given by $T_\ast \phi(\alpha_p)
= (T\phi^{-1})^\ast \alpha_p = \alpha_p \circ T_{\phi(p)} \phi^{-1}$
for $\alpha_p \in T_p^\ast M.$ We further write $\B \phi$ for the
inverse of $\F \phi$. With this notation, the following proposition
describes all automorphisms of the standard Courant algebroid, see
\cite[Prop.~3.24]{gualtieri:2003a}:
\begin{proposition}
    \label{theo.couralgauto}
    Let $E = \stCA$ be the standard Courant algebroid. Then every
    automorphism $\Phi$ of $E$ is of the form
    \begin{equation}
        \label{eq:Automorphisms}
        \Phi = \tau_B \circ
        \F\phi,
    \end{equation}
    with a unique closed $2$-form $B \in \Omega^2(M)$ and an unique
    diffeomorphism $\phi:M \ra M$. The automorphism group of $\stCA$
    is given by the semi-direct product $\mathcal{Z}^2(M) \rtimes
    \mathcal{D}\text{\itshape{iff}\/}(M)$ with $\mathcal{Z}^2(M) =
    \op{ker}(\dif_{|\Omega^2(M)})$, where the group multiplication is
    \begin{equation}
        \label{eq:SemiDirect}
        (B,\phi)(C,\psi) = (B + (\phi^{-1})^\ast C, \phi \circ \psi).
    \end{equation}
\end{proposition}

%
%

\subsection{Dirac Structures}
\label{subsec:DiracStructures}

The definition of Dirac structures on manifolds is due to Courant
\cite{courant:1990a} and was later generalized to Courant algebroids
in \cite{liu.weinstein.xu:1997a}:
\begin{definition}
    Let $E$ be a Courant algebroid. A subbundle $L \subset E$ is
    called a Dirac structure if $L$ is maximally isotropic with
    respect to the given bilinear form and if $\sect(L)$ is closed
    under the Courant bracket, i.e. $[\sect(L),\sect(L)] \subseteq
    \sect(L)$.
\end{definition}
In the following, whenever we speak about Courant algebroids with
Dirac structures, we will restrict ourself to Courant algebroids with
\emph{even fiber dimension} and a bilinear form of \emph{signature
  zero}.  The reason for this is that maximal isotropic subbundles in
such Courant algebroids have half the fiber dimension of the
algebroid, a point that will become important later on.

In particular, the standard Courant algebroid $\stCA$ is of this type.
In this case, one has the following two standard examples of Dirac
structures:
\begin{example}
    \label{example:PresymplecticPoisson}
    Let $E =\stCA$ be the standard Courant algebroid over $M$.
    \begin{enumerate}
    \item Given a two-form $\omega \in \Omega^2(M)$, we consider
        $\omega$ as a map $\omega: TM \ra T^\ast M$ be defining
        \begin{equation}
            \label{eq:OmegaOfX}
            \omega(X) = i_X \omega = \omega(X,\ldot).
        \end{equation}
        Thanks to skew-symmetry the $\dim M$-dimensional subbundle $L
        := \graph(\omega) \subset TM\oplus T^\ast M$ is isotropic.
        Moreover $L$ is closed under the Courant bracket, i.e.  a
        Dirac structure, if and only if $\omega$ is closed. Thus
        presymplectic two-forms can be viewed as particular cases of
        Dirac structures.
    \item Let $\pi \in \sect(\Wedge^2 TM)$ be a bivector. We consider
        $\pi$ as a map $\pi: T^*M\ra TM$ by defining
        \begin{equation}
            \label{eq:PiOfAlpha}
            \pi(\alpha) = \pi(\alpha,\ldot).
        \end{equation}
        Again due to skew-symmetry $L := \graph(\pi) \subset \stCA$ is
        a maximal isotropic subbundle. One further finds that $L$ is a
        Dirac structure if and only if $\pi$ is a Poisson tensor, i.e.
        $[\pi,\pi] = 0$.
    \end{enumerate}
\end{example}

%
%

\section{Derived Brackets for Courant Algebroids and Dirac Structures}
\label{sec:DerivedBrackets}

In this section we shall realize the Courant bracket as a derived
bracket in the sense of
\cite{kosmann-schwarzbach:1996a,kosmann-schwarzbach:2004b} as this has
been done before by Roytenberg \cite{roytenberg:1999a,
  roytenberg:2002b} in a slightly different context.

%
%

\subsection{The Rothstein-Poisson Bracket}
\label{subsec:RothsteinBracket}

For the study of Poisson manifolds the Schouten-Nijenhuis bracket has
turned out to be a very useful tool since one can write the Poisson
bracket as a derived bracket $\{f,g\} = -[[f,\pi],g]$ for a unique
bivector $\pi \in \sect(\Wedge^2 TM)$. It then follows immediately
that the Jacobi identity for the Poisson bracket is equivalent to the
equation $[\pi,\pi] = 0$, see e.g.\cite{kosmann-schwarzbach:2004b} for
an overview on derived brackets. In the case of a Courant algebroid
$E$ a similar approach is possible. However, one first has to find an
appropriate space which has the sections $\sect(E)$ as a subset as
well as a bracket on it, in order to write the Courant bracket as a
derived bracket. One possibility favored by Roytenberg
\cite{roytenberg:1999a, roytenberg:2002b} is given by the space of
functions on a suitable symplectic supermanifold.

We shall use a slightly different presentation avoiding the explicit
notion of supermanifolds: in our approach we take advantage of more
conventional differential geometry by using the Rothstein-Poisson
bracket \cite{rothstein:1991a} on the sections of the Grassmann
algebra of the pullback bundle $\tau^\varsharp E \longrightarrow
T^*M$, see Appendix~\ref{sec:Rothstein} for precise definitions.  The
Rothstein-Poisson bracket satisfies a graded Leibniz rule with respect
to the $\wedge$-product, is graded antisymmetric and fulfilles a
graded Jacobi identity where all signs come from the Grassmann parity.
Though the structure is essentially the same as in
\cite{roytenberg:1999a, roytenberg:2002b}, which can made even more
transparent in the super-Darboux coordinates from
Section~\ref{section.superdarboux}, the explicit use of ordinary
differential geometry might come in useful when considering a
quantized version of Dirac structures as we can rely on e.g.
Bordemann's construction\cite{bordemann:2000a, bordemann:1996a:pre}
for deformation quantization of the Rothstein-Poisson bracket.
Furthermore, the usage of the Rothstein-Poisson bracket allows us to
perform intrinsically \emph{global} computations.

Let $E \ra M$ be a vector bundle together with a fiber metric $h$,
i.e. a nondegenerate bilinear form, and let $\nabla$ be a metric
connection on $E$. We denote by $\tau: T^*M \longrightarrow M$ the
cotangent bundle. Then on the supercommutative algebra
$\grassalg{\tau^\varsharp E}$ of sections of the pulled back bundle
$\tau^\varsharp E \longrightarrow T^*M$ we have the Rothstein-Poisson
bracket as described in Appendix~\ref{sec:pullback_Rothstein}, defined
by use of the pulled back of the fiber metric $h$ and the connection
$\nabla$.  We can regard $\grassalg{E}$ as a subalgebra of
$\grassalg{\tau^\varsharp E}$ via the pull-back of sections.

Since $T^*M$ is a vector bundle itself and since we consider a pulled
back bundle over $T^*M$, it makes sense to speak of sections $e \in
\grassalg{\tau^\varsharp E}$ which are polynomial in the fiber
directions of $T^*M$ of degree $k \in \mathbb{N}$.  Note that the
grading with respect to the fiber variables (the momenta) is not a
good grading for the Rothstein-Poisson bracket, neither is the
Grassmann degree. However, the Rothstein-Poisson bracket is graded
with respect to twice the polynomial degree in the momenta plus the
Grassmann degree.  We denote homogeneous sections of this \emph{total
  degree} $k \in \mathbb{N}$ by $\pol^k \subseteq
\grassalg{\tau^\varsharp E}$.  Then their direct sum $\pol^\bullet$ is
a subalgebra of $\grassalg{\tau^\varsharp E}$, both with respect to
the $\wedge$-product and the Rothstein bracket. With respect to this
grading, the Rothstein-Poisson bracket has degree $-2$, i.e.
\begin{equation}
    \label{eq:RothsteinGrading}
    \roth{\pol^k, \pol^\ell} \subseteq \pol^{k+\ell-2}.
\end{equation}
In particular, $C^\infty(M) = \pol^0$ and $\tau^\varsharp\sect(E) =
\pol^1$, see Appendix~\ref{grading}.

%
%

\subsection{Courant Algebroids via Rothstein Bracket}
\label{subsec:CourantViaRothstein}

With the help of the Rothstein-Poisson bracket on
$\grassalg{\tau^\varsharp E}$ we can define a derived bracket
\cite{kosmann-schwarzbach:1996a} on $\sect(E)$.  Consider for $\Theta
\in \grassalg{\tau^\varsharp E}$ the bilinear map on
$\grassalg{\tau^\varsharp E}$ given by
\begin{equation}
    \label{eq:DerivedBracket}
    (\xi,\zeta)
    \longmapsto 
    \roth{\roth{\xi,\Theta},\zeta}.
\end{equation}
In order to get a derived bracket on $\sect(E)$, the subspace $\pol^1
= \sect(E)$ has to be closed under the above map.  As one can see in
the local formula \eqref{equation.rothstein_pullback_bundle} for the
Rothstein-Poisson bracket this is only the case for a homogeneous
$\Theta \in \pol^3$ of total degree $3$.  In fact, the lower degrees
do not contribute and higher ones will not produce pullbacks of
sections from $\sect(E)$. Such a section $\Theta \in \pol^3$ has two
types of contributions: one is a section of
$\Gamma^\infty(\tau^\varsharp E)$ which is \emph{linear} in the
momenta variables of $T^*M$, the other is a pull-back section of
$\Gamma^\infty (\bigwedge^3 E)$.
\begin{lemma}
    Let $\Theta \in \pol^3$. Then the following objects are well-defined:
    \begin{enumerate}
    \item A $\R$-bilinear derived bracket $\der{\Theta}{\dcd}:
        \sect(E) \times \sect(E) \ra \sect(E)$ defined for $e_1,e_2
        \in \sect(E)$ by
        \begin{equation}
            \label{eq:ThetaDerivedBracket}
            \der{\Theta}{e_1, e_2} =  \roth{\roth{e_1, \Theta}, e_2}.
        \end{equation}
    \item A derived anchor, i.e. a bundle map $\anc{\Theta}:E \ra TM$
        defined for $e\in\sect(E)$ and $f\in C^\infty(M)$ by
        \begin{equation}
            \label{eq:DerviedAnchor}
            \anc{\Theta}(e)f = \roth{\roth{e, \Theta}, f}.
        \end{equation}
    \item A map $\DT:C^\infty (M) \ra \sect(E)$ defined for $f\in
        C^\infty(M)$ by
        \begin{equation}
            \label{eq:DTDef}
            \DT f = \roth{\Theta, f}.
        \end{equation}
    \end{enumerate}
    The bundle $E$ together with the bilinear form $h$ and the above
    defined bracket and anchor satisfy the conditions {\it ii.)} and
    {\it iii.)} from definition \ref{def.courant-algebroid} of a
    Courant algebroid.
\end{lemma}
\begin{proof}
    The well-definedness follows from the grading properties. Then the
    verification of the conditions \textit{ii.)}  and \textit{iii.)}
    is a straightforward computation using the graded Jacobi identity of
    $\roth{\cdot, \cdot}$.
\end{proof}
Note that the definition of $\DT$ is consistent with
Definition~\ref{def.courant-algebroid}.  The following lemma is the
analogue of \cite{roytenberg:2002b} for the Rothstein-Poisson bracket
and follows the general ideas of derived brackets
\cite{kosmann-schwarzbach:2004b}.
\begin{lemma}
    \label{lemma.squarezero}
    Let $\Theta \in \pol^3 \subset \grassalg{(\tau^\varsharp E)}$ be
    homogeneous of degree $3$.  Then $E$ together with the bilinear
    form $h$, the bracket $\der{\Theta}{\dcd}$ and the anchor
    $\anc{\Theta}$ is a Courant algebroid if and only if
    \begin{equation}
        \label{eq:ThetaThetaNull}
        \roth{\Theta,\Theta} =0.   
    \end{equation}
\end{lemma}
\begin{proof}
    For the `if' part we only have to check the Jacobi identity for
    $[\cdot, \cdot]$, which is a simple computation in the
    framework of derived brackets \cite{kosmann-schwarzbach:1996a}. We
    only have to use the graded Jacobi identity of $\roth{\cdot, \cdot}$. For
    the `only if' part we assume the Jacobi identity for
    $\derT{\dcd}$. Then
    \[
    \roth{\roth{\roth{\roth{\Theta, \Theta}, e_1}, e_2}, e_3} = 0
    \tag{$*$}
    \]
    for all $e_1,e_2,e_3 \in \sect(E)$.  Let $f \in C^\infty(M)$ be a
    function. Then by the graded Leibniz rule for $\roth{\cdot, \cdot}$ we
    have
    \[
    0 = \roth{\roth{\roth{\roth{\Theta, \Theta}, e_1}, e_2}, f e_3}
    = f \roth{\roth{\roth{\roth{\Theta, \Theta}, e_1}, e_2}, e_3} +
    \roth{\roth{\roth{\roth{\Theta, \Theta}, e_1}, e_2}, f} e_3
    \]
    from which we obtain
    \[
    \roth{\roth{\roth{\roth{\Theta, \Theta}, e_1}, e_2}, f} = 0.
    \tag{$**$}
    \]
    By another application of the graded Jacobi identity we also find
    $\roth{\roth{\roth{\Theta, \Theta}, \roth{e_1, e_2}}, f} = 0$ for
    all $e_1,e_2 \in \sect(E)$ and $f \in C^\infty(M)$. Since locally
    every function $g \in C^\infty(M)$ can be written as $g =
    \roth{e_1, e_2}$ with appropriate $e_1, e_2 \in \sect(E)$ we
    conclude
    \[
    \roth{\roth{\roth{\Theta, \Theta}, f}, g} = 0
    \tag{$*$$**$}
    \]
    for all $f,g \in C^\infty(M)$. From the explicit formulas for
    $\roth{\cdot, \cdot}$ we see that the properties ($*$), ($**$) and
    ($*$$**$) together imply that the homogeneous element
    $\roth{\Theta, \Theta}$ of degree $4$ has to vanish.
\end{proof}

In a next step we want to construct such an element $\Theta$ for a
given Courant algebroid. We begin with the following easy lemma:
\begin{lemma}
    \label{lemma:TheTorsionTensor}
    Let $(E, \cour{\dcd}, \rho, h)$ be a Courant algebroid with a
    metric connection $\nabla$. Then the map $T:
    \sect(E)\times\sect(E)\times\sect(E) \ra \R$ defined by
    \begin{equation}
        \label{eq:TDef}
        T(e_1,e_2,e_3) = h(\nabla_{\rho(e_1)}e_2 - \nabla_{\rho(e_2)}
        e_1 - \cour{e_1,e_2},e_3) + h(\nabla_{\rho(e_3)} e_1,e_2)
    \end{equation}
    is a skew-symmetric $3$-tensor $T \in \sect(\Wedge^3 E^\ast)$.
\end{lemma}
\begin{proof}
    The proof is a direct computation using the definition of a
    Courant algebroid and the fact that the connection is metric.
\end{proof}

In some sense, $T$ is the Courant algebroid version of the torsion of
$\nabla$.  Let $\coorddown{u}{K}$ be a local basis of sections of $E$
with dual basis $u^1, \ldots, u^K$, defined on the domain of a local
chart $x^1, \ldots, x^n$ of $M$.  Then locally $T$ is given by
\begin{equation}
    \label{eq:Tlocally}
    T = \frac{1}{6} T_{ABC}u^A\wedge u^B \wedge u^C    
\end{equation}
where $T_{ABC} = T(u_A, u_B, u_C)$. Using the musical isomorphism
$\sharp$ induced by the fiber metric $h$ we obtain from $T \in
\sect(\Wedge^3 E^*)$ the tensor field $T^\sharp \in \sect(\Wedge^3
E)$, locally given by
\begin{equation}
    \label{eq:Tduallocally}
    T^\sharp 
    = \frac{1}{6} h^{AE} h^{BF} h^{CG} T_{ABC} u_E\wedge u_F \wedge u_G.    
\end{equation}
Using the structure functions $C_{AB}^C = \scal{\cour{u_A,u_B},u^C}$
of the Courant bracket, the components $\rho^i = \dif x^i\circ \rho$ of
the anchor, and the Christoffel symbols $\Gamma^A_{iB}$ of $\nabla$ we
obtain by a straightforward computation
\begin{equation}
    \label{eq:TdualLocal}
    T^\sharp = 
    \frac{1}{2} h^{AD} h^{BE} \rho^i(u_A) 
    \Gamma^C_{iB}  u_C\wedge u_D\wedge u_E 
    -
    \frac{1}{6}  h^{AD} h^{BE} C_{AB}^C u_C\wedge u_D \wedge u_E.
\end{equation}
The second tensor field we shall need is obtained as follows. Since
the anchor can be viewed as $\rho \in \Gamma^\infty(E^* \otimes TM)$
we can use $h$ to obtain a tensor field $\rho^\sharp \in
\Gamma^\infty(E \otimes TM)$. In a second step we can view the tangent
vector field part of $\rho^\sharp$ as a linear function on $T^*M$
whence we end up with a section $\mathcal{J}(\rho^\sharp) \in
\Gamma^\infty(\tau^\varsharp E)$, polynomial in the momenta of degree
$1$. Here, $\mathcal{J}: \Gamma^\infty(\mathrm{S}^\bullet TM)
\longrightarrow \Pol^\bullet(T^*M)$ denotes the canonical algebra
isomorphism. Locally, $\mathcal{J}(\rho^\sharp)$ is given by
\begin{equation}
    \label{eq:JofRhoSharp}
    \mathcal{J}(\rho^\sharp) = h^{AC} p_i \rho^i(u_A) u_C,
\end{equation}
where $p_1, \ldots, p_n$ are the canonically conjugate momenta on
$T^*M$ to the local coordinates $q^1 = \tau^*x^1, \ldots, q^n = \tau^*
x^n$ induced by the local coordinates $x^1, \ldots, x^n$ on $M$.

Putting both together we obtain from the choice of a metric connection
$\nabla$ the homogeneous element
\begin{equation}
    \label{eq:ThetaNabla}
    \Theta = - \mathcal{J}(\rho^\sharp) + T^\sharp \in \pol^3
    \subseteq \Gamma^\infty(\tau^\varsharp \Wedge^\bullet E)
\end{equation}
of total degree $3$. For later use we shall give yet another local
expression for $\Theta$, namely using the super-Darboux coordinates
from Proposition~\ref{proposition:SuperDarbouxCoordinates}. By
rearranging the local expressions for $\mathcal{J}(\rho^\sharp)$ and
$T^\sharp$ we obtain
\begin{equation}
    \label{eq:ThetaSuperDarboux}
    \Theta = 
    - h^{AC} r_i \rho^i(u_A) u_C 
    - \frac{1}{6} h^{AD} h^{BE} C^C_{AB} u_C\wedge u_D \wedge u_E.
\end{equation}
The advantage will be the easy commutation relations between the $r_i$
and the other local variables. It is also the direct analogue to the
supergeometric formulation of Roytenberg, see
\cite[Eq.~(4.7)]{roytenberg:2002b}. Note however, that this splitting
of $\Theta$ is \emph{not} coordinate independent, i.e. the two parts
are \emph{not} tensor fields, contrary to the splitting
\eqref{eq:ThetaNabla}.
\begin{lemma}
    \label{lemma:ThetaFromCourant}
    Let $E \ra M$ be a Courant algebroid and chose a metric connection
    $\nabla$. Define the element $\Theta \in \pol^3$ by
    \eqref{eq:ThetaNabla}. Then the Courant bracket and the anchor of
    $E$ coincide with the derived bracket and the derived anchor
    induced by the element $\Theta \in \pol^3$. In particular,
    $\roth{\Theta, \Theta} = 0$.
\end{lemma}
\begin{proof}
    Using the super-Darboux coordinates this is a simple verification.
    The second statement follows directly from
    Lemma~\ref{lemma.squarezero}.
\end{proof}

Now we can finally make contact to the supermanifold formulation of
Roytenberg. Analogously to \cite[Thm. 4.5]{roytenberg:2002b} we
obtain: 
\begin{theorem}
    \label{theorem.derived_bracket}
    Let $E \longrightarrow M$ be a vector bundle with fiber metric $h$
    and metric connection $\nabla$. Then the set of Courant algebroid
    structures on $E$ is in one-to-one correspondence with the set of
    $\Theta \in \pol^3$ such that $\roth{\Theta, \Theta} =0$.
\end{theorem}

%
%

\subsection{The Case $E = L\oplus L^\ast$}

Consider the case $E = L\oplus L^\ast$ for a vector bundle $L$ endowed
with the natural pairing as fiber metric of signature zero.  In the
following we shall use a connection on $L$ and the corresponding
induced metric connection on $L \oplus L^*$. From this choice we
obtain the Rothstein-Poisson bracket on $\spalg$, see also
Appendix~\ref{subsec:caseL+L*}. The splitting $E = L \oplus L^*$
induces a bigrading instead of our previous total degree: Indeed, we
set $\degL$ to be the polynomial degree in the momenta plus the
$L$-degree and $\degLs$ is the polynomial degree in the
momenta plus the $L^*$-degree. Then $\pol^{(r,s)}$ denotes those
elements in $\pol^{r+s}$ of $\degL$-degree $r$ and $\degLs$-degree
$s$. Using this direct sum decomposition one obtains the following,
analogously to \cite{roytenberg:2002a}:
\begin{lemma}
    \label{lemma.homological_long}
    Let $\Theta = \psi + \mu + \gamma +\phi \in \pol$ be an element of
    total degree $3$ with $\psi \in \pol^{(0,3)}$,
    $\mu\in\pol^{(1,2)}$, $\gamma \in \pol^{(2,1)}$ and
    $\phi\in\pol^{(3,0)}$. Then $\roth{\Theta,\Theta} = 0$ is
    equivalent to
    \begin{eqnarray*}
        \roth{\mu,\psi} &=& 0\\
        \frac{1}{2}\roth{\mu,\mu} + \roth{\gamma,\psi} &=& 0\\
        \roth{\phi,\psi} + \roth{\mu,\gamma} &=&0\\
        \frac{1}{2} \roth{\gamma,\gamma} + \roth{\mu,\phi}& =& 0\\
        \roth{\gamma,\phi} &=& 0.
    \end{eqnarray*}
\end{lemma}

A quasi-Lie algebroid is a vector bundle $A \ra M$ together with a
$\mathbb{R}$-bilinear, skew symmetric bracket $\bra{A}{\dcd}$ on
$\sect(A)$ and a vector bundle homomorphism $\anc{A}: A\ra TM$ such
that for all $a_1,a_2 \in \sect(A)$ and $f\in C^\infty(M)$ the Leibniz
rule
\begin{equation}
    \label{eq:LeibnizQuasi}
    \bra{A}{a_1,f a_2} 
    = f \bra{A}{a_1,a_2} + \anc{A}(a_1)f \; a_2
\end{equation}
is satisfied.  If in addition the Jacobi identity for the bracket
$\bra{A}{\dcd}$ is fulfilled, then $A$ is a Lie algebroid.  In this
case the anchor $\anc{A}$ is a homomorphism of Lie algebras,
\begin{equation}
    \label{eq:AnchorMorph}
    \anc{A}(\bra{A}{a_1,a_2}) = [\anc{A}(a_1),\anc{A}(a_2)].
\end{equation}
The construction of the Schouten-Nijenhuis bracket can be generalized
to the quasi-Lie algebroid case by imposing the graded Leibniz rule.
The graded skew-symmetry and Leibniz rule are still satisfied, the
graded Jacobi identity holds if and only if $A$ is a Lie algebroid.
The Lie algebroid differential also generalizes to the quasi-Lie
algebroid case, and we get a differential with square zero if and only
if $A$ is a Lie algebroid. Conversely the quasi-Lie algebroid
structure can be retrieved from the Schouten-Nijenhuis bracket or from
the differential, see e.g.~\cite{marle:2002a}. Identifying our derived
brackets in this situation gives the following explicit formulas:
\begin{lemma}
    \label{lemma.derived_lie_brackets}
    Let $\Theta = \psi + \mu + \gamma +\phi \in \pol \subset \spalg$
    be as above, and let $\der{\Theta}{\dcd}$ and $\anc{\Theta}$ be
    the derived bracket and anchor.
    \begin{enumerate}
    \item The restriction of $\der{\Theta}{\dcd}$ to $L$ with subsequent
        projection to $L$ is given by the derived bracket with respect
        to $\mu$,
        i.e. for all $s_1,s_2 \in \sect(L)$ we have
        \begin{equation}
            \label{eq:BracketForL}
            \op{pr}_L(\der{\Theta}{s_1,s_2}) 
            = \der{\mu}{s_1,s_2} 
            = \roth{\roth{s_1,\mu},s_2}.
        \end{equation}
        Further, the restriction of the
        anchor to $L$ is given by 
        \begin{equation}
            \label{eq:AnchorForL}
            \anc{\Theta}(s) f = \anc{\mu}(s) f =
            \roth{\roth{s,\mu},f}
        \end{equation}
        for $s \in \sect(L)$ and $f \in C^\infty(M)$.

    \item The bracket $\der{\mu}{\dcd}$ together with the anchor
        $\anc{\mu}$ make $L$ a quasi Lie-algebroid. The associated
        Schouten-Nijenhuis bracket is given by
        \begin{equation}
            \label{eq:SchoutenBracket}
            \der{\mu}{P,Q}=\roth{\roth{P,\mu},Q}
        \end{equation}
        for $P,Q \in \sect(\Wedge^\bullet L)$, and the Lie algebroid
        differential by
        \begin{equation}
            \label{eq:DifferentialForL}
            \dif_L \eta = \roth{\mu,\eta},
        \end{equation}
        where $\eta \in \sect(\Wedge^\bullet L^\ast).$
    \item Analogous results are obtained for $L^\ast$  by replacing
        $\mu$ with $\gamma$.
    \end{enumerate}
\end{lemma}
\begin{proof}
    Let $s_1,s_2 \in \pol^{(1,0)} = \sect(L)$. Using the bigrading
    properties we get
    \[
    \der{\Theta}{s_1,s_2} = \roth{\roth{s_1,\Theta},s_2}
    = \roth{\roth{s_1,\psi},s_2} + \roth{\roth{s_1,\mu},s_2}
    \]
    with $\roth{\roth{s_1, \psi}, s_2} \in \pol^{(0,1)} =
    \sect(L^\ast)$ and $ \roth{\roth{s_1, \mu}, s_2} \in \pol^{(1,0)}
    = \sect(L)$. Thus $\op{pr}_L(\der{\Theta}{s_1, s_2}) =
    \der{\mu}{s_1, s_2} = \roth{\roth{s_1, \mu}, s_2}$. Analogously,
    we obtain \eqref{eq:AnchorForL}. A standard computation finally
    shows that the extension of $\der{\mu}{\cdot, \cdot}$ to
    multivector fields is given by \eqref{eq:SchoutenBracket} since
    \eqref{eq:SchoutenBracket} satisfies the same type of graded Leibniz rule
    and coincides with the Schouten-Nijenhuis bracket on the local
    generators.  As one can see by counting degrees we have a
    well-defined map $\roth{\mu, \ldot}: \sect(\Wedge^k L^\ast) \ra
    \sect(\Wedge^{k+1} L^\ast)$. Thanks to the graded Leibniz rule for the
    Rothstein-Poisson bracket, this map is a graded derivation of the
    $\wedge$-product. A straightforward computation then shows that
    $i_s \roth{\mu, f} = i_s \dif_L f$ and $i_{s_2} i_{s_1}
    \roth{\mu,\alpha} = i_{s_2} i_{s_1} \dif_L \alpha$ for all $s,
    s_1, s_2 \in \sect(L)$, $f \in C^\infty (M)$ and $\alpha \in
    \sect(L^\ast)$.  By the derivation property, $\roth{\mu, \cdot}$
    coincides with $\dif_L$ on the whole space $\sect(\Wedge^\bullet
    L^\ast)$. The last statement follows analogously.
\end{proof}

Recall that a Lie quasi-bialgebroid is a Lie algebroid
$(A,\bra{A}{\dcd},\anc{A})$ together with a graded derivation
$\dif_{A^*}$ of degree one of $\grassalg{A}$ with respect to both the
$\wedge$-product and the Schouten-Nijenhuis bracket, and a $3$-vector
$\phi \in \sect(\Wedge^3 A)$ such that $\dif_{A^*} \phi = 0$ and
$\dif_{A^*}^2 = -\bra{A}{\phi,\cdot}$, see
e.g.~\cite{laurent-gengoux.ponte.xu:2005a}.  It is well-known that a
graded derivation of degree one of $\grassalg{A}$ defines a
Schouten-Nijenhuis bracket on $\grassalg{A^\ast}$. Thus a Lie
quasi-bialgebroid is a pair $(A,A^\ast)$, where $A$ is a Lie algebroid
and $A^\ast$ is a quasi-Lie algebroid, such that the differential
$\dif_{A^\ast}$ of $A^\ast$ is a graded derivation of the
Schouten-Nijenhuis bracket on $A$ and $\dif_{A^\ast}^2 =
-\bra{A}{\phi,\cdot}$ for some $3$-vector $\phi \in \sect(\Wedge^3 A)$
with $\dif_{A^\ast} \phi = 0$, see e.g.~\cite{roytenberg:2002a}.  A
Lie bialgebroid \cite{mackenzie.xu:1994a} is obtained in the case that
$\dif_{A^\ast}^2=0$.  Combining Lemma~\ref{lemma.homological_long} and
\ref{lemma.derived_lie_brackets} we obtain the following well-known
result \cite{roytenberg:1999a}:
\begin{lemma}\label{lemma.quasiLie}
    Let $\Theta = \psi + \mu + \gamma +\phi \in \pol \subset \spalg$
    be as in Lemma \ref{lemma.derived_lie_brackets} and assume now in
    addition $\roth{\Theta,\Theta} = 0$. If $\psi = 0$ then $L$ is a
    Lie quasi-bialgebroid. If $\phi = 0$ then $L^\ast$ is a Lie
    quasi-bialgebroid.
\end{lemma}

%
%

\subsection{Courant Algebroids with Dirac Structures}
\label{subsec:CourantAlgsWithDiracStructures}

We shall now consider the case of a Courant algebroid $E = L\oplus
L^\ast$ over $M$ such that $L$ is a Dirac structure. As we will see
later in Lemma \ref{corollary_E=L+J(L)} a Courant algebroid $E$ with a
Dirac structure $L$ is always of this form.

The element $\Theta$ from Theorem~\ref{theorem.derived_bracket} now is
given as a sum $\Theta = \psi +\mu+\gamma+\phi$ according to the
bigrading.  We split the tensor field $T$ from
Lemma~\ref{lemma:TheTorsionTensor} and $T^\sharp$, respectively, into
their $L$ and $L^*$ components.  Note also that we can identify $T$
and $T^\sharp$ canonically, since $L \oplus L^*$ is canonically
`self-dual'.  As before we set $\rho^i = \dif x^i \circ \rho$ and
define $r_i \in \sect(\Wedge^\bullet \tau^\sharp (L\oplus L^\ast))$
for $i=1, \ldots, n$ by $r_i = p_i - \Gamma^\beta_{i \alpha}
a^\alpha\wedge a_\beta$, see
Proposition~\ref{proposition:SuperDarbouxCoordinatesL+L*}.
Analogously, the anchor $\rho$ splits into the two restrictions
$\rho_L$ and $\rho_{L^*}$ to $L$ and $L^*$, respectively. Then we have
locally
\begin{equation}
    \label{eq:rhoLLocally}
    \mathcal{J}(\rho_L^\sharp) = p_i \rho^i(a_\alpha) a^\alpha
    \quad
    \textrm{and}
    \quad
    \mathcal{J}(\rho_{L^*}^\sharp) = p_i \rho^i(a^\alpha) a_\alpha
\end{equation}
The above splitting of $T$ and $\rho$ into the components according to
$E = L \oplus L^*$ now gives the splitting of $\Theta$ into the
elements $\mu$, $\gamma$, and $\phi$. To identify these components, we
\emph{define} the global tensor fields
\begin{equation}
    \label{eq:Defmu}
    \mu = 
    - \mathcal{J}(\rho_L^\sharp) 
    + T\big|_{\bigwedge^2 L \otimes L^*}
\end{equation}
\begin{equation}
    \label{eq:Defgamma}
    \gamma = 
    - \mathcal{J}(\rho_{L^*}^\sharp) 
    + T\big|_{L \otimes \bigwedge^2 L^*}
\end{equation}
\begin{equation}
    \label{eq:phiDef}
    \phi = T\big|_{\bigwedge^3 L^\ast}.
\end{equation}
Because $L$ is a Dirac structure one has $T\big|_{\bigwedge^3 L}=0$
and therefor
\begin{equation}
    \Theta = \mu +\gamma + \phi.
\end{equation}

A little computation shows that $ T\big|_{\bigwedge^2 L \otimes L^*}$
is three times the torsion \cite{crainic.fernandes:2003a} for the Lie
algebroid $L$ and analogously $T\big|_{L \otimes \bigwedge^2 L^*}$ is
three times the torsion for the quasi Lie algebroid $L^\ast$. We
further have
\begin{equation}
    \phi(\sigma_1,\sigma_2,\sigma_3) =
    -\scal{\cour{\sigma_1,\sigma_2},\sigma_3}
\end{equation}
for $\sigma_1$, $\sigma_2$, $\sigma_3 \in \pol^{(0,1)} =
\sect(L^\ast)$.

Let us look now at the local expressions. Let $\coordup{x}{n}$ be
coordinates on $M$, $\coorddown{a}{k}$ a local basis of sections of
$L$ and $\coordup{a}{k}$ the dual local basis of sections of $L^\ast$.
We define local functions
\begin{equation}
    \label{eq:StructureFunktionDirac}
    c_{\alpha\beta}^\gamma 
    = \scal{\cour{a_\alpha, a_\beta},a^\gamma}
    \quad
    \textrm{and}
    \quad
    c^{\alpha\beta}_\gamma 
    = \scal{\cour{a^\alpha,a^\beta},a_\gamma}.
\end{equation}
Furthermore, we have
\begin{equation}
    \label{eq:phiLocalExpression}
    \phi = \frac{1}{6}\phi^{\alpha\beta\gamma} 
    a_\alpha\wedge a_\beta \wedge a_\gamma
    \quad
    \textrm{with}
    \quad 
    \phi^{\alpha\beta\gamma} 
    = -(\scal{\cour{a^\alpha,a^\beta},a^\gamma}).
\end{equation}
Since we assume $L$ to be a Dirac structure, all other combinations of
structure functions as in Section~\ref{subsec:CourantViaRothstein} are
either zero, or can be computed from the ones in
\eqref{eq:StructureFunktionDirac} and \eqref{eq:phiLocalExpression} by
using the properties of the Courant bracket $\cour{\cdot, \cdot}$.

Let $\coordup{q}{n},\coorddown{p}{n}$ be the induced coordinates on
$T^\ast M$ and let $\Gamma^\alpha_{i \beta}$ be the Christoffel
symbols for the connection $\nabla$ on $L$.  We define
\begin{equation}
    \label{eq:TorsionT}
    T_{\alpha\beta}^\gamma 
    = T(a_\alpha, a_\beta, a^\gamma)
    = \rho^i(a_\alpha) \Gamma_{i \beta}^\gamma 
    - \rho(a_\beta)^i \Gamma_{i \alpha}^\gamma 
    - c_{\alpha\beta}^\gamma
\end{equation}
\begin{equation}
    \label{eq:TorsionTbar}
    T^{\alpha\beta}_\gamma 
    = T(a^\alpha, a^\beta, a_\gamma)
    = \rho^i(a^\beta) \Gamma^\alpha_{i\gamma}
    - \rho^i(a^\alpha) \Gamma_{i \gamma}^\beta 
    - c^{\alpha\beta}_\gamma
\end{equation}
and then we have locally
\begin{equation}
    T = \frac{1}{2} T_{\alpha \beta}^\gamma a^\alpha\wedge a^\beta\wedge a_\gamma
    + \frac{1}{2} T^{\alpha \beta}_\gamma a_\alpha\wedge a_\beta\wedge
    a^\gamma + \frac{1}{6} \phi^{\alpha\beta\gamma}  a_\alpha\wedge
    a_\beta\wedge a_\gamma.
\end{equation}
>From this we immediately obtain the following statement:
\begin{lemma}
    Locally, $\mu$ and $\gamma$ are given by
    \begin{equation}
        \label{eq:MuLocally}
        \mu = -p_i \rho^i(a_\alpha)a^\alpha 
        + \frac{1}{2} T_{\alpha\beta}^\gamma 
        a^\alpha\wedge a^\beta\wedge a_\gamma
    \end{equation}
    \begin{equation}
        \label{eq:GammaLocally}
        \gamma = -p_i \rho^i(a^\alpha)a_\alpha 
        + \frac{1}{2} T^{\alpha\beta}_\gamma
        a_\alpha\wedge a_\beta\wedge a^\gamma.
    \end{equation}
    In particular, $\mu \in \pol^{(1,2)}$ and $\gamma \in
    \pol^{(2,1)}$. In the local super-Darboux coordinates we have
    \begin{equation}
        \label{eq:MuDarboux}
        \mu = - r_i \rho^i(a_\alpha) a^\alpha 
        - \frac{1}{2} c_{\alpha\beta}^\gamma 
        a^\alpha\wedge a^\beta\wedge a_\gamma
    \end{equation}
    \begin{equation}
        \label{eq:GammaDarboux}
        \gamma = -r_i \rho^i(a^\alpha) a_\alpha 
        - \frac{1}{2} c^{\alpha\beta}_\gamma 
        a_\alpha\wedge a_\beta\wedge a^\gamma.
    \end{equation}
\end{lemma}
Using the splitting in \eqref{eq:MuDarboux} and
\eqref{eq:GammaDarboux} we can compare with Roytenberg's expressions
in \cite[Eq.~(3.10) and (3.11)]{roytenberg:1999a}. Note however, that
this splitting depends on the choice of coordinates while
\eqref{eq:MuLocally} and \eqref{eq:GammaLocally} have intrinsic
geometric meanings.
\begin{corollary}
    \label{Theta}
    We have $\roth{\Theta,\Theta} = 0$, or equivalent
    \begin{eqnarray*}
        \roth{\mu,\mu} &=& 0\\
        \frac{1}{2} \roth{\gamma,\gamma} + \roth{\mu,\phi}& =& 0\\
        \roth{\mu,\gamma} &=&0\\
        \roth{\gamma,\phi} &=& 0.
    \end{eqnarray*}
\end{corollary}
\begin{example}
    Let $E=\stCA$ be the standard Courant algebroid over $M$. Let
    $\nabla$ be any torsion-free connection and construct the
    Rothstein-Poisson bracket on $\sect(\Wedge^\bullet\tau^\varsharp
    E)$.  First we get
    \begin{equation}
        \label{eq:TMTstarMalleNull}
        \gamma = 0,
        \quad
        \phi = 0,
        \quad
        \textrm{and}
        \quad
        \psi = 0.
    \end{equation}
    For the only nontrivial element $\mu$ we find locally $\mu = - p_i
    \tau^\varsharp \dif x^i \in \sect(\tau^\varsharp T^*M)$. The
    pulled back bundle $\tau^\varsharp T^*M$ can be identified with
    the \emph{annihilator subbundle}
    $\mathrm{Ver}(T^*M)^{\mathrm{ann}} \subseteq T^*(T^*M)$ of the
    \emph{vertical} subbundle $\mathrm{Ver}(T^*M) \subseteq T(T^*M)$
    in the usual way. This canonical identification allows us to
    identify $\tau^\varsharp \dif x^i$ with $\tau^* \dif x^i = \dif
    q^i$. Hence, under this identification, $\mu$ coincides with the
    \emph{canonical one-form} $- \theta_0$ on $T^*M$.
\end{example}
We have some more corollaries to Corollary~\ref{Theta} and
Lemma~\ref{lemma.quasiLie}. First we note \cite{roytenberg:1999a}:
\begin{corollary}
    \label{corollary:DiracQuasiLieBiAlg}
    On $L$ we have given the structure of a quasi-Lie bialgebroid.
\end{corollary}
\begin{remark}
    If in addition $\roth{\mu, \phi} = 0$ is satisfied then $(L,
    L^\ast)$ is a Lie bialgebroid \cite{kosmann-schwarzbach:2005a}.
    But only if $\phi = 0$ the space of sections $\sect(L^\ast)$ is
    closed under the Courant bracket and $L^\ast$ is a Dirac
    structure.
\end{remark}
Given a Dirac structure $L$ in a Courant algebroid $E$ we always can
find a maximal isotropic subbundle $L'$ complementary to $L$ and
identify $E$ with $L\oplus L^\ast$, see e.g.
Corollary~\ref{corollary_E=L+J(L)}. Thus we have
\cite{roytenberg:1999a}:
\begin{corollary}
    \label{corollary:EisLoplusLstar}
    A Courant algebroid $E$ with a Dirac structure $L$ is isomorphic
    to the double of the Lie quasi-bialgebroid $L\oplus L^\ast$.
\end{corollary}

As shown in Lemma~\ref{lemma.derived_lie_brackets} the derived bracket
$\der{\mu}{\dcd}$ is the Schouten-Nijenhuis bracket for the Lie
algebroid structure on $L$ given by the restriction of the Courant
bracket and the anchor to $\sect(L)$. Further $\der{\gamma}{\dcd}$
defines a quasi-Lie algebroid structure on $L^\ast$ where the bracket
is given by
\begin{equation}
    \label{eq:BracketLstar}
    \bra{L^\ast}{\sigma_1, \sigma_2} 
    = \der{\gamma}{\sigma_1, \sigma_2} 
    = \op{pr}_{L^\ast}(\cour{\sigma_1, \sigma_2})
\end{equation}
and the anchor by $\anc{L^\ast} = \rho|_{L^\ast}$.  The differential
$\dif_L$ is a graded derivation for the bracket $\der{\gamma}{\dcd}$
and the differential $\dif_{L^\ast}$ is a graded derivation for the
bracket $\der{\mu}{\dcd}$.

%
%

\section{Smooth and Formal Deformations of Dirac Structures}
\label{sec:SmoothFormalDeformations}

In this section we shall now establish the smooth and the formal
deformation theory of Dirac structures. In the following $E$ is a
Courant algebroid with a fiber metric of signature zero and $L
\subseteq E$ a Dirac structure as before.

%
%

\subsection{Definition of Smooth Deformations}
\label{subsec:DefinitionSmooth}

As motivation we first recall the well-known situation for Poisson
manifolds, see e.g. \cite[Sect.~18.5]{cannasdasilva.weinstein:1999a}:
A \emph{smooth deformation} $\pi_t$ of a Poisson structure $\pi_0$ on
$M$ is a smooth map
\begin{equation}
    \label{eq:SmoothDefPoisson}
    \pi: I \times M \ra \Wedge^2 TM
\end{equation}
with $\pi_t = \pi(t,\ldot) \in \sect(\Wedge^2 TM)$ for all $t\in
I$ and $\pi(0,\ldot)=\pi_0$, such that 
\begin{equation}
    \label{eq:SmoothDefPoissonJacobi}
    [\pi_t,\pi_t] = 0    
\end{equation}
for all $t \in I$, where $I\subseteq \R$ is an open interval around
zero.  \emph{Formal} deformations then are given by formal power
series $\pi_t = \pi_0 + t \pi_1 + \ldots \in \sect(\Wedge^2 TM)[[t]]$
such that $[\pi_t,\pi_t] = 0$ order by order in the formal parameter.
A similar approach is possible in the case of symplectic manifolds.

Consider now a Dirac structure $L$ in $E$.  One possibility to define
a \emph{smooth} deformation of $L$ is given by specifying a family of
subbundles in terms of a family of projections. This way, we can
encode the desired smoothness easily:
\begin{definition}
    Let $L \subseteq E$ be a Dirac structure and let $I \subseteq \R$
    be an open interval around zero. A smooth deformation of $L = L_0$
    is a family of Dirac structures $L_t$ with $t \in I$ such that
    there exists a smooth map
    \begin{equation}
        \label{eq:Pdef}
        P:I\times M \ra \End(E)
    \end{equation}
    with
    \begin{enumerate}
    \item $P(t,m):  E_m \ra E_m $ for all $t \in I$ and $m \in M$
    \item $P(t,m)^2 = P(t,m)$ for all $t \in I$ and $m \in M$
    \item $\op{Im} P_t = L_t$ for all $t \in I$, where $P_t =
        P(t,\ldot) \in \sect(\End(E))$. 
    \end{enumerate}
\end{definition}
\begin{remark}
    \label{remark:AlternativeToSmoothDef}
    Consider the pull-back bundle $\op{pr}^\varsharp E$, where
    $\op{pr}: I \times M \ra M$ is the projection. Equivalent to the
    definition above we can consider a smooth deformation of $L$ as a
    smooth subbundle $\mathfrak{L} \subseteq \op{pr}^\varsharp E$ such
    that every $L_t = \mathfrak{L}|_{\{t\}\times M} \subset E$ is a
    Dirac structure where $L_0 = L$.
\end{remark}
While the above definition is conceptually clear and easy, it is not
very suited for concrete computations. Thus we shall re-formulate the
definition using additional geometric structures in
Section~\ref{subsec:RewritingDefProblem}. We also have to discuss the
possible notions of equivalence in detail. However, we first recall
two general well-known properties of the subbundles in question:
\begin{theorem}
    \label{theorem.fastcomplex}
    Let $E$ be a vector bundle with a fiber metric $\bili{\dcd}$. Then
    there exits a positive definite fiber metric $g$ and an isometry
    $J: E \ra E$ of $\bili{\dcd}$ with $J^2 = \id$, such that
    \begin{equation}
        \label{eq:NichtKomplexeStruktur}
        g(e_1,e_2) = (e_1,J e_2)
    \end{equation}
    for all $e_1,e_2 \in \sect(E)$.
\end{theorem}
\begin{proof}
    For the readers convenience we sketch the proof: Choose a positive
    definite fiber metric $k$ and define $A \in
    \Gamma^\infty(\End(E))$ by $k(A e_1,e_2) = \bili{e_1,e_2}$.  Since
    $A$ turns out to be $k$-symmetric we can use its polar
    decomposition $A = \sqrt{A^2} J$. Then $g(e_1,e_2) = \bili{e_1,J
      e_2}$ has the required properties.
\end{proof}
\begin{corollary}
    \label{corollary_E=L+J(L)}
    Let $E$ be a vector bundle with even fiber dimension $2k$ and let
    $\bili{\dcd}$ be a bilinear form on $E$ of signature zero. Let further
    be $L$ a maximal isotropic subbundle of $E$. Choose $g$ and $J$
    according to Theorem~\ref{theorem.fastcomplex}. Then
    \begin{equation}
        \label{eq:ELoplusJL}
        E = L\oplus J(L)
        \quad
        \textrm{and}
        \quad
        L^{\perp_g} = J(L) \cong L^\ast.
    \end{equation}
\end{corollary}
\begin{theorem}
    \label{theorem.fedosov}
    Let $E$ be a vector bundle, $I \subset \R$ an open interval around
    zero and let $L_t$ for $t\in I$ be a smooth family of subbundles
    of $E$. Then there exits a vector bundle automorphism $U_t$ of $E$
    over the identity $\id: M \longrightarrow M$, 
    smoothly depending on $t\in I$ such that 
    \begin{equation}
        \label{eq:Ut}
        L_t = U_t(L_0).
    \end{equation}
    If $E$ is a Courant algebroid and $L_t$ a family of maximal
    isotropic subbundles, then we can also achieve that $U_t$ is an
    isometry of the symmetric bilinear form $h = \bili{\dcd}$ for all
    $t \in I$.
\end{theorem}
\begin{proof}
    The theorem can be proved along the lines of
    \cite[Lem.~1.1.5]{fedosov:1996a}.
\end{proof}

%
%

\subsection{The Problem of Equivalence}
\label{subsec:ProblemEquivalence}

Let $L$ be a Dirac structure in a Courant algebroid $E$ and $L_t$ a
smooth deformation of $L=L_0$. We know from
Theorem~\ref{theorem.fedosov} that there exists an isometry $U_t$ of
$E$ smoothly depending on $t$ such that $L_t = U_t(L_0)$. Thus, in
this general concept it seems natural to define a \emph{trivial
  deformation} as a deformation $L_t$ such that we can find a time
dependent $U_t$ which is not only an isometry but also a \emph{Courant
  algebroid automorphism}. If we further ask whether two smooth
deformations $L_t$ and $L'_t$ are equivalent, one is tempted to
require the existence of a time dependent Courant algebroid
automorphism $U_t$, such that $L'_t = U_t(L_t)$.

However, in the case of of the standard Courant algebroid $\stCA$, due
to Lemma~\ref{theo.couralgauto}, this would mean that we have the
gauge transformations by closed two-forms as equivalence
transformations. But then every two Dirac structures given by
presymplectic forms would be equivalent. Hence we see that in the case
of $\stCA$ we can not permit every Courant algebroid automorphism as
equivalence transformation as long as we want to reproduce the common
results for the deformation theory of symplectic forms.

In the case of $E = \stCA$, we know that every automorphism is given
by the the product of a gauge transformation and a lifted
diffeomorphism $\F \phi$. As we do not want gauge transformations as
equivalence transformations we have to consider the \emph{lifted
  diffeomorphisms}.  Indeed, given a presymplectic form $\omega$ on a
manifold $M$ and a diffeomorphism $\phi$ of $M$, one can easily show
\cite{bursztyn.radko:2003a} that the equation
\begin{equation}
    \label{eq:Backphiomega}
    \B\phi(\graph \omega) = \graph(\phi^\ast \omega)
\end{equation}
is satisfied.  Analogously we have 
\begin{equation}
    \label{eq:Backphipi}
    \B\phi(\graph \pi) = \graph(\phi^\ast \pi)    
\end{equation}
for a Poisson tensor $\pi$ on $M$. This motivates the following
definition of equivalent deformations of Dirac structures which reduce
to the well-known situation in the Poisson or symplectic case:
\begin{definition}
    \label{definition:Equivalence}
    Let $L \subset \stCA$ be a Dirac structure in the standard Courant
    algebroid. Two smooth deformations $L_t$ and $L'_t$ of $L$ are
    called equivalent, if there exists a smooth curve of
    diffeomorphisms $\phi_t$ of $M$ such that $L'_t = \F \phi_t(L_t)$.
    A smooth deformation is called trivial, if there exists a smooth
    curve of diffeomorphisms $\phi_t$ such that $L_t = \F\phi_t(L_0)$.
\end{definition}

While for the standard Courant algebroid this seems to be the
reasonable definition of equivalent deformations, in general it will
be more difficult: for any vector bundle $E \longrightarrow M$ we have
the exact sequence of groups
\begin{equation}
    \label{eq:Sequence}
    1 \longrightarrow \mathrm{Gau}(E) \longrightarrow \mathrm{Aut}(E) 
    \longrightarrow \mathrm{Diffeo}(M) \longrightarrow 1,
\end{equation}
where $\mathrm{Gau}(E)$ denotes those vector bundle automorphisms of
$E$ which induce the identity on $M$, and the last arrow assigns to an
arbitrary vector bundle automorphism $\Phi: E \longrightarrow E$ the
induced diffeomorphism $\phi$ of $M$. However, quite unlike for the
Courant algebroid $TM \oplus T^*M$, in general this exact sequence
does \emph{not} split. Furthermore, even if the sequence splits, it is
not clear, whether the split can be chosen in a reasonable way.  In
fact, if $E$ is associated to the frame-bundle, then one can choose a
splitting.

Since it is precisely this canonical splitting in the case of $TM
\oplus T^*M$ which we use for Definition~\ref{definition:Equivalence}
there seems to be no simple way out. One possibility would be the
following: since we are only interested in smooth curves of
diffeomorphisms $\phi_t$ of $M$ with $\phi_0 = \id_M$ we know that
such a curve is the time evolution of a time-dependent vector field
$X_t$ on $M$. After the \emph{choice} of a connection $\nabla$ on $E$
we can lift $X_t$ horizontally to $E$ and consider its time evolution
$\Phi_t$ on $E$. Then we can use $\Phi_t$ instead of $\F \phi_t$ to
formulate a definition of equivalence and trivial deformations
analogously to Definition~\ref{definition:Equivalence}. However, this
would depend explicitly on the choice of a connection. We shall come
back to this problem in a future work.  At the present stage, the
question of equivalence of smooth deformations of Dirac structures in
a general Courant algebroid has to be left unanswered.

%
%

\subsection{Rewriting the Deformation Problem}
\label{subsec:RewritingDefProblem}

To study the formal deformation theory of Dirac structures we first
have to think about an appropriate description for such deformations.
Given a Courant algebroid $E$ with Dirac structure $L$ we
\emph{choose} an isotropic complement $L'$ to $L$ (for example with
the help of Corollary~\ref{corollary_E=L+J(L)}) and identify $L'$ with
$L^\ast$.  Then we can write $E = L \oplus L^\ast$, where the fiber
metric on $E$ translates to the natural pairing on $L \oplus L^\ast$.
Thus we may assume that $E$ has this form in the following. Note
however, that we still have to discuss the influence of this chosen
isomorphism later.

Locally a small deformation $L_t$ of $L$ could be understood as the
graph of a map $\omega_t: L \ra L^\ast$. Indeed, over a compact subset
$K \subseteq M$ a smooth deformation $L_t$ can be written as the graph
of some $\omega_t$ provided $t$ is sufficiently small. Globally in
$M$, this needs not to be true whence smooth deformation theory
becomes highly non-trivial. However, since we will mainly be
interested in formal deformations (to be thought of as formal Taylor
expansions of smooth deformations) the idea of looking at graphs will
be sufficient for us.  The claim that $L_t$ is isotropic allows us to
identify $\omega_t$ with a $2$-form in $L$.  To ensure that
$\sect(L_t)$ is closed under the Courant bracket and therefor is a
Dirac structure leads to an additional requirement for $\omega_t$.

In the following consideration we will first omit the dependency on
$t$. So let $\omega \in \Omega^2(L)$ be a $2$-form. Then
$\op{graph}(\omega)$ is integrable, i.e. closed under the Courant
bracket, if and only if for all $s_1, s_2, s_3 \in \sect(L)$
\begin{equation}
    \label{equation.deformation_equation.long}
    \begin{split}
        0 
        &= 
        \scal{
          \cour{s_1 + \omega(s_1), s_2 + \omega(s_2)},
          s_3 + \omega(s_3)} \\
        &=
        \scal{\cour{s_1, \omega(s_2)}, s_3} 
        + \scal{\cour{\omega(s_1), s_2}, s_3} 
        + \scal{\cour{s_1, s_2}, \omega(s_3)}\\
        &\quad 
        + \scal{\cour{s_1, \omega(s_2)}, \omega(s_3)}
        + \scal{\cour{\omega(s_1), s_2]}, \omega(s_3)} 
        + \scal{\cour{\omega(s_1), \omega(s_2)}, s_3} \\
        &\quad
        + \scal{\cour{\omega(s_1), \omega(s_2)}, \omega(s_3)}.
    \end{split}
\end{equation}
The constant term in $\omega$ vanishes as $L$ is assumed to be a Dirac
structure throughout. Moreover, this equation combines linear,
quadratic and cubic terms in $\omega$.  In order to analyze this
equation in more detail, we use the Rothstein-Poisson bracket.
\begin{lemma} 
    \label{lemma:DeformatioEquation}
    Let $E = L \oplus L^\ast$ be a Courant algebroid with $L$ a Dirac
    structure and let $\omega \in \sect(\Wedge^2 L^\ast)$ be a
    $2$-form. Then $\op{graph}(\omega)\subseteq E$ is a Dirac
    structure if and only if
    \begin{eqnarray}
        \label{equation.deformation_equation.roth}
        \roth{\mu, \omega} 
        + \frac{1}{2} \roth{\roth{\omega, \gamma}, \omega} 
        + \frac{1}{6} 
        \roth{\roth{\roth{\phi, \omega}, \omega}, \omega} 
        = 0. 
    \end{eqnarray}
\end{lemma}
\begin{proof}
    Replace all Courant brackets by the derived bracket using
    $\Theta$ gives \eqref{equation.deformation_equation.roth} after a
    straightforward computation.
\end{proof}

Due to the bigrading properties of the Rothstein-Poisson bracket the
definition
\begin{equation}
    \label{eq:TripleBracket}
    [\eta_1,\eta_2,\eta_3]_\phi =
    \roth{\roth{\roth{\phi, \eta_1},\eta_2},\eta_3}
\end{equation}
gives a well-defined trilinear map
\begin{equation}
    \label{eq:TripleBracketMap}
    \sect(\Wedge^k L^\ast) 
    \times \sect(\Wedge^l L^\ast)
    \times \sect(\Wedge^m L^\ast) 
    \longrightarrow
    \sect(\Wedge^{k+l+m-3} L^\ast).
\end{equation}
Moreover, because $\phi$ is a pull-back section this map is
independent of the connection used for constructing the
Rothstein-Poisson bracket. With the definitions from
Lemma~\ref{lemma.derived_lie_brackets} we can write
\eqref{equation.deformation_equation.roth} equivalently as
\begin{equation}
    \label{equation.deformation_equation.short}
    \dL \omega 
    + \frac{1}{2}\der{\gamma}{\omega,\omega} 
    + \frac{1}{6}[\omega,\omega,\omega]_\phi =
    0.
\end{equation}
This is the fundamental equation for $\omega$ which has been derived
by Roytenberg in his approach in another context, see
\cite{roytenberg:2002a}.

Equation \eqref{equation.deformation_equation.short} is precisely the
sorting of \eqref{equation.deformation_equation.long} by the
homogeneous monomials in $\omega$ and hence independent of the usage
of the Rothstein-Poisson bracket. Nevertheless, we can use the
Rothstein-Poisson bracket to obtain algebraic identities for the three
parts of \eqref{equation.deformation_equation.short} which are very
hard to obtain without the Rothstein-Poisson bracket.

%
%

\subsection{Formal Deformations}
\label{subsec:FormalDeformations}

Following the general idea of formal deformation theory, namely to
solve a non-linear algebraic equation order by order in terms of
formal power series \cite{gerstenhaber.schack:1988a,
  gerstenhaber:1964a}, we consider solutions of
\eqref{equation.deformation_equation.long} in the sense of formal
power series. Since $\omega$ should be a `small' deformation we make
the Ansatz
\begin{equation}
    \label{eq:omegaAnsatz}
    \omega = t \omega_1 + t^2 \omega_2 + \cdots = \sum_{t=1}^\infty
    t^r \omega_r \in t \Gamma^\infty(\Wedge^2 L^*)[[t]],
\end{equation}
where $\omega_1, \omega_2, \ldots$ have to be determined recursively.
Since $\omega_r \in \Gamma^\infty(\Wedge^2 L^*)$, we can interprete
the deformation as a $2$-cochain in the Lie algebroid complex of $L$,
viewed only as a Lie algebroid.  The following lemma is now crucial
for the cohomological approach:
\begin{lemma}
    \label{lemma.universalid}
    Let $\eta \in \sect(\Lambda^2 L^\ast)$ be a two-form. Then
    \begin{equation}
        \dif_\eta 
        = \roth{\mu,\ldot} 
        + \roth{\roth{\eta, \gamma}, \cdot} 
        + \frac{1}{2}\roth{\roth{\roth{\phi, \eta}, \eta}, \cdot}
        =
        \dif_L 
        + \der{\gamma}{\eta,\cdot} + \frac{1}{2}[\eta,\eta,\cdot]_\phi
    \end{equation}
    defines a graded derivation of degree one of the
    $\wedge$-product such that
    \begin{equation}
        \label{eq:CooleGleichung}
        \dif_\eta \left(
            \dif_L \eta
            + \frac{1}{2}\der{\gamma}{\eta, \eta}
            + \frac{1}{6}[\eta,\eta,\eta]_\phi
        \right) = 0.
    \end{equation}
\end{lemma}    
\begin{proof}
    Using the derived bracket formalism this is a straightforward
    computation.
\end{proof}

The following theorem shows that the solvability of
\eqref{equation.deformation_equation.long} or equivalently
\eqref{equation.deformation_equation.short} order by order in the
formal parameter leads to a cohomological obstruction in the usual
way:
\begin{theorem}
    Let $E = L \oplus L^\ast$ be a Courant algebroid with a Dirac
    structure $L$ and let $\omega_t = t \omega_1 + t^2\omega_2 +
    \ldots + t^N \omega_N \in \sect(\Wedge^2 L^\ast)[[t]] $ be a
    formal deformation of $L$ of order $N$, i.e. the equation
    \begin{equation}
        \label{eq:OrderN}
        \dL \omega_t 
        + \frac{1}{2}\der{\gamma}{\omega_t, \omega_t}
        +\frac{1}{6}[\omega_t, \omega_t, \omega_t]_\phi= 0
    \end{equation}
    is satisfied up to order  $N$. Then
    \begin{equation}
        \label{eq:Rest}
        R_{N+1} =
        - \frac{1}{2}
        \sum_{i=1}^{N} \der{\gamma}{\omega_i, \omega_{N+1-i}} 
        - \frac{1}{6}
        \sum_{i+j+k=N+1} [\omega_i, \omega_j, \omega_k]_\phi
        \in \Gamma^\infty(\Wedge^3 L^*)
    \end{equation}
    is closed with respect to $\dL$, and $\omega_t$ can be extended to
    a deformation of order $N+1$ if and only if $R_{N+1}$ is exact.
\end{theorem}
\begin{proof}
    The proof is essentially the usual argument of formal deformation
    theory.  Let $\omega_{N+1} \in \sect(\Wedge^2 L^\ast)$ be
    arbitrary and set $\omega'_t = \omega_t + t^{N+1}\omega_{N+1}$.
    Then 
    \[
    \begin{split}
        \dL \omega'_t 
        & + \frac{1}{2}\der{\gamma}{\omega'_t, \omega'_t}
        + \frac{1}{6} [\omega'_t, \omega'_t, \omega'_t]_\phi \\ 
        &=
        t^{N+1}\left(
            \dL \omega_{N+1} 
            + \frac{1}{2}\sum_{i=1}^{N} 
            \der{\gamma}{\omega_i, \omega_{N+1-i}} 
            +  \frac{1}{6}\sum_{i+j+k=N+1}
            [\omega_i, \omega_j, \omega_k]_\phi 
        \right) + o(t^{N+2}),
    \end{split}
    \]
    whence $\omega'_t$ satisfies
    \eqref{equation.deformation_equation.short} up to order $N+1$ if
    and only if $\dL \omega_{N+1} = R_{N+1}$, i.e. $R_{N+1}$ is exact
    with respect to $\dL$. On the other hand, $R_{N+1}$ is always
    closed. Indeed, by Lemma~\ref{lemma.universalid} applied to
    $\omega_t'$ we get
    \[
    0 = t^{N+1}\dL
    \left(
        \frac{1}{2} \sum_{i=1}^{l}
        \der{\gamma}{\omega_i, \omega_{N+1-i}} 
        + \frac{1}{6} \sum_{i+j+k=N+1}
        [\omega_i, \omega_j, \omega_k]_\phi
    \right)
    + o(t^{N+2}),
    \]
    which implies $\dL R_{N+1} = 0$.
\end{proof} 
\begin{remark}
    \label{remark:WhyNontrivial}
    From the proof it is clear that the whole derived bracket
    formalism enters only in showing that $R_{N+1}$ is $\dL$-closed.
    This is in some sense the nontrivial statement of the theorem. In
    principle, this can also be shown directly using only
    \eqref{equation.deformation_equation.long} and the algebraic
    identities for the Courant bracket. However, the computations are
    very much involved without using the nice derived bracket
    formalism. Nevertheless, it should be emphasized that the
    characterization of the order-by-order obstruction to solve
    \eqref{equation.deformation_equation.long} by the third Lie
    algebroid cohomology of the Dirac structure is \emph{independent}
    of the choices we made in order to obtain the Rothstein-Poisson
    bracket.
\end{remark}

%
%

\subsection{Examples: Presymplectic and Poisson Manifolds}
\label{subsec:Examples}

Let us now discuss some examples in order to show that the deformation
theory of Dirac structures generalizes the well-known deformation
theories of presymplectic and Poisson structures.

Let $(M,\omega)$ be a presymplectic manifold, and consider the
standard Courant algebroid $\stCA$ with the Dirac structure $L =
\graph(\omega)$. In this case $T^\ast M$ is a complement of $L$ and we
can identify $L^\ast \cong T^\ast M$. There is also a canonical
identification of $L$ with $TM$, which is given by the restriction of
the gauge transformation $\tau_{-\omega}$ to $L$.  Because $\omega$ is
closed, $\tau_{-\omega}$ is a Courant algebroid automorphism and we
have the identification $L\oplus L^\ast \cong \stCA$, where the
Courant algebroid structure on the right hand side is still the
standard one. Because $L^\ast \cong T^\ast M$ is a Dirac structure
with trivial Lie algebroid structure, according to our theory a smooth
deformation of $L\cong TM$ is given by a closed time-dependent
two-form $\eta_t$ with $\eta_0 = 0$. The deformation of the original
Dirac structure is then given by $L_t = \op{graph}(\omega +\eta_t)$,
i.e. by the deformation of the presymplectic form $\omega$. Thus we
retrieve the common results in this case.

Usually in formal deformation theory, the infinitesimally inequivalent
deformations are parameterized by a second cohomology relevant for the
deformation problem while the third cohomology gives the obstructions
for the existence of order-by-order deformations. In our case, one
would expect the second Lie algebroid cohomology to be the relevant
one.

For the usual deformation theory of symplectic forms or Poisson
bivectors this is indeed the case. However, in the general case, the
situation is more subtle. To see this, we consider the following
example:

First recall that the Lie algebroid cohomology of a Dirac structure
coming from a presymplectic structure coincides with the de Rham
cohomology. Then, for a presymplectic manifold, two formal
deformations $\omega_t$ and $\omega'_t$ of the presymplectic form
$\omega_0$ are equivalent iff there exists a formal diffeomorphism
$\phi_t = \op{exp}(\Lie_{X_t})$ with $X_t = t X_1 + \ldots \in
t\sect(TM)[[t]]$, such that $\phi_t \omega_t = \omega'_t$. This is the
reasonable definition of `deformations up to formal diffeomorphisms'.
In first order this equation reads as
\begin{equation}
    \label{eq:OmegasEquivalent}
    \omega'_1 - \omega_1 = \Lie_{X_1}\omega_0 = \dif i_{X_1} \omega_0.
\end{equation}
If there is a $\alpha \in \Omega^1(M)$, such that $\dif\alpha =
\omega'_1 - \omega_1$, then we must find $X_1$ with $i_{X_1}\omega_0 =
\alpha$. For a \emph{symplectic} form $\omega_0$ this is always
possible, so nontrivial deformations only exist if $H^2_{\op{dR}}(M)$
is nontrivial. However, if we start with a presymplectic form
$\omega_0$, there might be \emph{no} $X_1$ such that $i_{X_1}\omega_0 =
\alpha$ and the triviality of $H^2_{\op{dR}}(M)$ is not sufficient for
the rigidity of $M$ as a presymplectic manifold. Because the
presymplectic deformation is a special case of the deformation of
Dirac structures, the obstructions for the existence of non-trivial
deformations are not in the second Lie algebroid cohomology of $L$.
This is probably the most surprising feature of the deformation theory
of Dirac structures.
\begin{remark}
    \label{remark:OtherOption}
    One might wonder whether this is just an artifact of our notion of
    equivalence based on formal diffeomorphisms. However, if one
    decides to use the notion of equivalence suggested by
    Theorem~\ref{theorem.fedosov} (which we do not prefer, see the
    discussion in Section~\ref{subsec:ProblemEquivalence}), then the
    situation is even worse: All deformations of presymplectic forms
    in this sense become equivalent, while the second Lie algebroid
    cohomology might be nontrivial.
\end{remark}

Finally, let us consider a Poisson manifold $(M,\pi)$, and consider
the Dirac structure $L = \graph(\pi)$ in the standard Courant
algebroid $\stCA$. We choose $TM$ as complement to $L$ so that we can
identify $L^\ast$ with $TM$.  Observe that in this case $L^\ast \cong
TM$ is again a Dirac structure but unlike as above the Lie algebroid
structure on $L^\ast$ is non-trivial. We further identify $L$ with
$T^\ast M$ via $\rho^\ast|_L$. Hence, we have the identification
$L\oplus L^\ast = T^\ast M\oplus TM$, but the Courant algebroid
structure on the right hand side now is not the standard one.  The
differential $\dif_L$ becomes the differential given by $\pi$, i.e.
$\dif_{\pi} = [\pi,\ldot]$, and the bracket on $L^\ast \cong TM$ is
the canonical Schouten-Nijenhuis bracket.  Deformations of $L$ are
given by time-dependent bivector fields $\lambda_t$ such that
\begin{equation}
    \label{eq:DeformationPoisson}
    \dif_{\pi} \lambda_t + \frac{1}{2}[\lambda_t,\lambda_t] = 0.
\end{equation}
Thanks to $[\pi,\pi] = 0$, this equation is equivalent to 
\begin{equation}
    \label{eq:DeformationPoissonII}
    [\pi + \lambda_t, \pi + \lambda_t]=0.
\end{equation}
We conclude that deformations $\pi_t = \pi + \lambda_t$ of the
Poisson tensor $\pi$ are the same as deformations of the corresponding
Dirac structure $L$.

%
%

\appendix

%
%

\section{The Rothstein-Poisson bracket}
\label{sec:Rothstein}

%
%

\subsection{Definition of the Rothstein-Poisson bracket}
\label{rothstein.general_definition}

Let $F \ra N$ be a vector bundle over a symplectic manifold
$(N,\omega)$ and let $\pi = -\omega^{-1}$ be the Poisson tensor for
the induced Poisson bracket on $N$, i.e. $\{f,g\} = \pi(\dif f,\dif
g)$. Further let $h$ be a pseudo-riemannian metric on $N$ and $\nabla$
a metric connection. We denote local coordinates on $N$ by
$\coordup{x}{n}$, local basis sections of $F$ by $\coorddown{s}{k}$
and the dual sections of $F^\ast$ by $\coordup{s}{k}$.  With the local
expression
\begin{equation}
    \hat{R} = \frac{1}{2}  \pi^{ij} h^{AB}
    R^C_{Ajk} \partial_i \otimes s_B\wedge s_C\otimes \dif x^k
\end{equation}
we get a well defined global section $\hat{R} \in \sect(TN \otimes
\Wedge^2 F \otimes T^\ast N)$, where $h^{AB}$ and $R^A_{Bij}$ are the
local expressions for the pseudo-riemannian metric $h^{-1}$ and the
curvature $R$ in coordinates.  A section $S \in \sect(TN \otimes
\Wedge^k F \otimes T^\ast N)$ can be interpreted as a map
\begin{equation}
    S:\sect(TN \otimes \Wedge^\bullet F) \ra \sect(TN \otimes
    \Wedge^{\bullet+k} F)
\end{equation}
by
\begin{equation}
    (X\otimes \phi  \otimes \eta)(Y\otimes \psi) 
    = \eta(Y) X\otimes \phi \wedge \psi.
\end{equation}
We therefor can form powers of $\hat{R}$ by composition of maps.
Because $\hat{R}$ increases the degree of the part in $\Wedge^\bullet
F$ by two, $\hat{R}$ is nilpotent and we have a well-defined section
\begin{equation}
    (\op{id} -\hat{R})^{-\frac{1}{2}} 
    = \op{id} + \frac{1}{2} \hat{R} +
    \frac{3}{8} \hat{R}^2 + \ldots,
\end{equation}
where $\op{id} = \partial_i \otimes 1 \otimes \dif x^i$ is the
identity map in $\sect(TN \otimes \Wedge^\bullet F)$.  For a section
$S \in \sect(TN \otimes \Wedge^\bullet F \otimes T^\ast N)$ we define
the local section $S^i_j$ of $\Wedge^\bullet F$ by
\begin{equation}
    S = \partial_i \otimes S^i_j \otimes \dif x^j.
\end{equation}
In the following we denote by $i(\sigma)\psi$ and $j(\sigma)\psi$ the
interior product of a section $\sigma \in \sect(F^\ast)$ with an
element $\psi \in\sect(\Wedge ^\bullet F)$ from the left and right,
respectivly.
\begin{theorem}[Rothstein-Poisson bracket \cite{rothstein:1991a, bordemann:2000a}]
    There is a super-Poisson bracket on $\grassalg{F}$ (depending on
    $\nabla,h$ and $\omega$) called the Rothstein-Poisson bracket
    which is locally given by
    \begin{equation}
        \roth{\phi,\psi} =
        \pi^{ij}\big(1-\hat{R})^{-\frac{1}{2}}\big)^k_i
        \wedge
        \big(1-\hat{R})^{-\frac{1}{2}}\big)^l_j 
        \wedge \nabla_{\partial_k} \phi 
        \wedge \nabla_{\partial_l} \psi 
        + h_{AB} j(s^A)\phi \wedge i(s^B)\psi.
    \end{equation}
\end{theorem}
That $\roth{\dcd}$ is a super-Poisson bracket means that for all $\phi
\in \sect(\Wedge^k F)$, $\psi \in \sect( \Wedge^l F)$ and $\eta \in
\sect( \Wedge^\bullet F) $ we have
\begin{enumerate}
\item $\roth{\phi,\psi} = -(-1)^{kl}\roth{\psi,\phi}$
\item $\roth{\phi,\psi\wedge \eta} = \roth{\phi,\psi}\wedge \eta +
    (-1)^{kl} \psi\wedge \roth{\phi,\eta}$
\item $\roth{\phi,\roth{\psi,\eta}} = \roth{\roth{\phi,\psi},\eta} +
    (-1)^{kl} \roth{\psi,\roth{\phi,\eta}}.$
\end{enumerate}

%
%

\subsection{The Rothstein bracket for pullback bundles}
\label{sec:pullback_Rothstein}

Let $\pi_M : E\ra M $ be a vector bundle, $N$ a manifold and let
$f:N\ra M$ be a smooth map. We denote by $f^\varsharp \pi_M:f^\varsharp E
\ra N$ the pullback bundle with respect to $f$. Given a section $e \in
\sect(E)$ the pullback section $f^\varsharp e \in \sect(f^\varsharp E)$ is
then defined by 
\begin{equation}
    f^\varsharp e = e\circ f.
\end{equation}
A local basis $u_1,\ldots,u_K \in \sect(E|_U)$ of $E$ defined on some
open set $U \subseteq M$ leeds to a local basis $f^\varsharp
u_1,\ldots,f^\varsharp u_K$ of $f^\varsharp E$ defined on
$f^{-1}(U)\subseteq N$.  Given a connection $\nabla$ on $E$ we have
the induced connection $f^\varsharp \nabla$ on $f^\varsharp E$, where
for pullback sections $f^\varsharp e$ with $e\in \sect(E)$ we have for
$Y \in \sect(TN)$
\begin{equation}
    f^\varsharp \nabla _Y(f^\varsharp e) 
    = f^\varsharp (\nabla_{Tf(Y)} e).
\end{equation}

In the following we will look at the case $N=T^\ast M$ with
$f=\tau:T^\ast M \ra M$ the cotangent projection. Let $x^1,\ldots,x^n$
be coordinates on $U\subseteq M$ and $q^1,\ldots,q^n,p_1,\ldots,p_n$
the induced bundle coordinates on $T^\ast U$.  The sign of the
canonical Poisson bracket on $T^\ast M$ is choosen such that
$\{q^i,p_j\} = \delta^i_j$.  Observe that for pullback sections $\tu
\in \tau^\varsharp(\sect(E))$ we have
\begin{equation}
    (\tau^\varsharp \nabla)_{\partdif{q^i}} \tu 
    = \tau^\varsharp (\nabla_{\partdif{x^i}} u)
    \quad
    \textrm{and}
    \quad
    (\tau^\varsharp \nabla)_{\partdif{p_i}} \tu = 0.
\end{equation}
In particular, for a general section $s \in \sect(\tau^\varsharp E)$
the expression $(\tau^\varsharp \nabla)_{\partdif{p_i}} s$ is
independent of the connection $\nabla$ and therefor we set
\begin{equation}
    \frac{\partial s}{\partial p_i} 
    = (\tau^\varsharp \nabla)_\partdif{p_i} s
\end{equation}
for the covariant derivative of $s$ with respect to $\partdif{p_i}$.
\begin{lemma}
    Let $\pi:E\ra M$ be a vector bundle with connection $\nabla^E$ and
    a fiber metric $h$. Look at the pullback bundle $F =
    \tau^\varsharp E$ over the symplectic manifold $T^\ast M$ together
    with the pullback connection $\nabla^F = \tau^\varsharp \nabla^E$
    and the pullback metric $\tau^\ast h$. Then the map $\hat{R}^F$ as
    defined in \eqref{rothstein.general_definition} satisfies
    \begin{eqnarray*}
        \hat{R}^F\Big(\partdif{q^i} \otimes \psi\Big) 
        &= &-\frac{1}{2}\partdif{p_j}\otimes \tau^\varsharp
        \( h^{AB}\: \(R^E\)^C_{Aij}\:u_B \wedge u_C\) 
        \wedge \psi\\ 
        \hat{R}^F\Big(\partdif{p_i} \otimes \psi\Big) &= &0,
    \end{eqnarray*}
    where $(R^E)^B_{Aij}$ is the curvature of $\nabla^E$ with respect to the
    appropriate coordinates.
\end{lemma}
\begin{proof}
    This is a straightforward computation.
\end{proof}
With this lemma it follows immediately that $(\hat{R}^F)^k = 0$ for
$k\geq 2$ and therefor
\begin{equation}
    (\op{id} - \hat{R}^F)^{-1} 
    = \op{id} + \frac{1}{2} \hat{R^F}.
\end{equation}
Hence in this case we get a more explicite formula for the
Rothstein-Poisson bracket.
\begin{lemma}
    With the above definitions the Rothstein-Poisson bracket on
    $\sect(\tau^\varsharp(\Wedge^\bullet E))$ is given by
    \begin{eqnarray}
        \label{equation.rothstein_pullback_bundle}
        \notag
        \roth{\phi,\psi} &=& \nabla^F_{\partdif{q^i}}\phi \wedge
        \partdif{p_i}\psi  - \partdif{p_i} \phi
        \wedge \nabla^F_{\partdif{q^i}} \psi  - \frac{1}{2}\tau^\varsharp\(
        h^{AB}\: \(R^E\)^C_{A ij}\: u_B \wedge u_C\) \wedge
        \partdif{p_i}\phi\wedge\partdif{p_j}\psi \\[2mm]
        && +\; \tau^\ast h_{AB}\; j(\tu^A)\phi \wedge i(\tu^B)\psi.
    \end{eqnarray}
\end{lemma}

%
%

\subsection{Super-Darboux coordinates}
\label{section.superdarboux}

Let us choose a local basis of sections $s_1,\ldots,s_k$ of the bundle
$E$ such that the functions $h_{AB}=h(s_A,s_B)$ are constant.  If we
calculate the Rothstein-Poisson bracket for the coordinate functions
$q^i,p_j$ and the local sections $\tu_A$ of $\tau^\varsharp E$ we get
the equations
\begin{equation}
    \begin{array}{rclrcl}
        \roth{q^i,q^j} &=& 0 & \roth{q^i,p_j} &=&
        \delta^i_j \\[2mm]
        \roth{p_i,p_j}&  =&-\frac{1}{2} \tau^\varsharp\( h^{AB}\: \(R^E\)^C_{A ij}\:
        u_B\wedge u_C\) \quad&
        \roth{q^i,\tu_A}& =& 0 \\[2mm]
        \roth{p_i,\tu_A}& =& -\tau^\varsharp\(\Gamma_{i
          A}^B\, u_B\) &  \roth{\tu_A,\tu_B}& =& \tau^\ast h_{AB}.
    \end{array}
\end{equation}
We see that $C^\infty(T^\ast M)$ is in general not closed under the
Rothstein-Poisson bracket. 
\begin{proposition}
    \label{proposition:SuperDarbouxCoordinates}
    Let the local sections $r_i$ of the bundle
    $\grassalg{(\tau^\varsharp E)}$ be defined by
    \begin{equation}
        r_i = p_i - \frac{1}{2} \tau^\varsharp\( h^{AB}\, \Gamma_{i A}^C\,
        u_B \wedge u_C\).
    \end{equation}
    Then the following equations are satisfied:
    \begin{equation}
        \roth{q^i,r_j} = \delta^i_j 
        \qquad \text{and} \qquad
        \roth{\tu_A,\tu_B} = h_{AB},
    \end{equation}
    and
    \begin{equation}
        \roth{q^i,q^j} = \roth{q^i,\tu_A} 
        =\roth{r_i,r_j} =  \roth{r_i,\tu_A} = 0.
    \end{equation}
\end{proposition}
\begin{proof}
    A direct calculation using the fact that the connection is metric
    leads to the result.
\end{proof}

%
%

\subsection{Grading for polynomial sections}
\label{grading}

Let $\pol \subset \grassalg{(\tau^\varsharp E)}$ be the sections which
are polynomial in the momenta, i.e. sections which locally can be
written as a linear combination of local sections of the form
\begin{equation}
    h^{A_1\ldots A_s} \tu_{A_1}\wedge\ldots\wedge \tu_{A_s}
\end{equation}
for $0\leq s \leq k$ with $h_{A_1\ldots A_s} \in \Pol^\bullet(T^\ast
M)$ polynomial functions on $T^\ast M$. From
\eqref{equation.rothstein_pullback_bundle} we get the following:
\begin{lemma}
    The space $\pol$ is closed under the Rothstein-Poisson bracket.
\end{lemma}
\begin{definition}
    Let the map $\deg:\pol \ra \pol$ be defined by the local formula
    \begin{equation}
        \deg = 2 p_i \partdif{p_i} + \tu_A \wedge
        i(\tu^A).
    \end{equation} 
    For an element $\phi \in \pol$ we say that $\phi$ is of degree $r$
    if the equation $\deg \phi = r \phi$ is satisfied. We denote the
    set of all such elements by $\pol^r$.
\end{definition}
\begin{remark}
    \begin{enumerate}
    \item Elements with degree zero can be identified with functions
        on $M$ and elements with degree one with sections in $E$, i.e.
        \begin{equation}
            \pol^0 = \tau^\ast(C^\infty(M))
            \quad
            \textrm{and}
            \quad
            \pol^1 = \tau^\varsharp(\sect(E)).
        \end{equation}
    \item The degree given by $\deg$ can be used to calculate the
        signs for the super-Poisson structure given by the Rothstein
        bracket because the momenta always count twice.
    \end{enumerate}
\end{remark}
\begin{lemma}
    The Rothstein-Poisson bracket is of degree $-2$ for the grading
    given by $\deg$.
\end{lemma}

%
%
 
\subsection{The case $E = L\oplus L^*$}
\label{subsec:caseL+L*}

Let $L \ra M$ be a vector bundle with a connection $\nabla$. We also
have a connection on the dual bundle $L^\ast$ and therefor a
connection $\nabla^E$ on $E = L \oplus L^*$, which is metric with
respect to the canonical bilinear form on $L\oplus L^*$, given by
\begin{equation}
    \scal{(s_1,\alpha_1),(s_2,\alpha_2)} = \alpha_1(s_2) +
    \alpha_2(s_1)
\end{equation}
for $s_1,s_2\in \sect(L)$ and $\alpha_1,\alpha_2 \in \sect(L^\ast)$.

Let $\coordup{x}{n}$ be coordinates on $M$, $a_1,\ldots,a_k$ be
a local basis of $L$ and $a^1,\ldots,a^k$ be the dual basis of
$L^\ast$.  Let $R_{\alpha i j}^\beta$ be the curvature on $L$ in
coordinates. The curvature on $L^\ast$ then is given in the dual
coordinates by $-R_{\alpha i j}^\beta$. If we choose
\begin{equation}
    (u_1,\ldots,u_A,\ldots,u_{2k}) = (a_1,\ldots
    a_k,a^1,\ldots a^k)
\end{equation}
as a local basis of $L \oplus L^\ast$ we get for the curvature on
$E=L\oplus L^\ast$
\begin{equation}
    \(R^E\)^B_{Aij} = \left\{
        \begin{array}{c@{\qquad}l}
            R^B_{Aij} &  \text{for}\; 1\leq A,B \leq k \\[2mm]
            -R^{A-k}_{B-k,ij} &\text{for}\; k+1\leq A,B \leq 2 k \\[2mm]
            0 & \text{otherwise}.\\
        \end{array} \right.
\end{equation}
Now let $F=\tau^\varsharp(L \oplus L^*) \longrightarrow T^*M$ again be
the pullback bundle.  Because of the special form of the curvature and
the fiber metric in the given coordinates, we can simplify
the formula for the Rothstein-Poisson bracket and get the following
lemma.
\begin{lemma}[Eilks \cite{eilks:2004a}]
    The Rothstein-Poisson bracket on $\spalg$ is locally given by
    \begin{eqnarray}
        \label{equation.rothstein_bracket.L+L*}
        \roth{\phi,\psi} &=& \nabla_{\partdif{q^i}}\phi \wedge
        \partdif{p_i}\psi  - \partdif{p_i} \phi
        \wedge \nabla_{\partdif{q^i}} \psi  + \tau^\varsharp\( 
        R^\alpha_{\beta ij} a_\alpha \wedge a^\beta\)\wedge
        \partdif{p_i}\phi\wedge\partdif{p_j}\psi \\[2mm] \notag
        && +\; j(\ta_\alpha)\phi \wedge i(\ta^\alpha)\psi +
        j(\ta^\alpha)\phi \wedge i(\ta_\alpha)\psi,
    \end{eqnarray}
    where $\psi$, $\phi \in \spalg$ and $\ta_\alpha$ and $\ta^\alpha$
    are pullback basis sections.
\end{lemma}
>From this formula we easily get the following lemma.
\begin{lemma}
    For sections $s \in \sect(L)$, $\sigma \in \sect(L^\ast)$,
     $P \in \sect(\Wedge^r L)$ and $\eta \in
    \sect(\Wedge^s L^\ast)$ we have the equations
    \begin{eqnarray*}
        \roth{\ts,\tau^\varsharp\eta} &=& \tau^\varsharp(i_s \eta) \\
        \roth{\tau^\varsharp \sigma,\tau^\varsharp P} &=& \tau^\varsharp(i_\sigma P)
    \end{eqnarray*}
    and
    \begin{equation}
        \roth{\ts,\tau^\varsharp P} 
        = 0 = 
        \roth{\tau^\varsharp \sigma,\tau^\varsharp \eta}.
    \end{equation}
    In particular, we get
    \begin{equation}
        \roth{\tau^\varsharp e_1,\tau^\varsharp e_2} =
        \tau^\ast\scal{e_1,e_2}
    \end{equation}
    for all $e_1,e_2 \in \sect(L\oplus L^\ast)$.
\end{lemma}
In this situation, the super-Darboux coordinates are given as follows:
\begin{proposition}
    \label{proposition:SuperDarbouxCoordinatesL+L*}
    If we set
    \begin{equation}
        r_i = 
        p_i - 
        \tau^\varsharp\left(
            \Gamma_{i\alpha}^\beta\, a^\alpha \wedge a_\beta\right),
    \end{equation}
    the only non-trivial Rothstein-Poisson brackets between the $q^i$,
    $r_j$, $\tau^\varsharp a^\alpha$ and $\tau^\varsharp a_\beta$ are
    \begin{equation}
        \roth{q^i,r_j} = \delta^i_j 
        \qquad 
        \text{and} 
        \qquad
        \roth{\ta^\alpha,\ta_\beta} = \delta^\alpha_\beta.
    \end{equation}
\end{proposition}
The grading with respect to the total degree can be refined in the
following sense:
\begin{definition}
    Let $\degL$ and $\degLs$ be defined by the local formula
    \begin{equation}
        \degL = p_i \partdif{p_i} + \ta_\alpha \wedge
        i(\ta^\alpha)
        \quad
        \textrm{and}
        \quad
        \degLs = p_i \partdif{p_i} + \ta^\alpha \wedge
        i(\ta_\alpha).
    \end{equation}
    For an element $\psi \in \pol$ we say $\psi$ has bidegree $(r,s)$,
    if $\degL \psi = r \psi$ and $\degLs\psi = s \psi$. The set of all
    such elements will be denoted by $\pol^{(r,s)}$.
\end{definition}
Of course we have $\deg = \degL + \degLs$, and therefor we call $\deg$
the total degree.  Moreover we have $\pol^{(0,0)} =
\tau^\varsharp(C^\infty(M))$, $\pol^{(r,0)} = \tau^\varsharp(\Wedge^r
L)$ and $\pol^{(0,s)} = \tau^\varsharp(\Wedge^s L^\ast)$.
\begin{lemma}
    The Rothstein-Poisson bracket restricted to the polynomial
    sections $\pol$ is of bidegree $(-1,-1)$, i.e. for all $\phi \in
    \pol^{(r,s)}$, $\psi \in \pol^{(t,u)}$ we have $\roth{\phi,\psi}
    \in \pol^{(r+t-1,s+u-1)}$.
\end{lemma}

%
%

\begin{footnotesize}
    \renewcommand{\arraystretch}{0.5}

\begin{thebibliography}{10}

\bibitem {bayen.et.al:1978a}
{\sc Bayen, F., Flato, M., Fr{{\o}}nsdal, C., Lichnerowicz, A., Sternheimer,
  D.: }\newblock {\em Deformation Theory and Quantization}.
\newblock Ann. Phys.  {\bf 111} (1978), 61--151.

\bibitem {bordemann:1996a:pre}
{\sc Bordemann, M.: }\newblock {\em On the deformation quantization of
  super-Poisson brackets}.
\newblock Preprint (Freiburg FR-THEP-96/8)  {\bf q-alg/9605038} (May 1996).

\bibitem {bordemann:2000a}
{\sc Bordemann, M.: }\newblock {\em The deformation quantization of certain
  super-{P}oisson brackets and {BRST} cohomology}.
\newblock In: {\sc Dito, G., Sternheimer, D. (eds.): }\newblock {\em
  Conf{\'e}rence Mosh{\'e} Flato 1999. Quantization, Deformations, and
  Symmetries}, {\em Mathematical Physics Studies} no. {\bf 22},   45--68.
  Kluwer Academic Publishers, Dordrecht, Boston, London, 2000.

\bibitem {bursztyn.radko:2003a}
{\sc Bursztyn, H., Radko, O.: }\newblock {\em Gauge Equivalence of {D}irac
  Structures and Symplectic Groupoids}.
\newblock Ann. Inst. Fourier  {\bf 53} (2003), 309--337.

\bibitem{cannasdasilva.weinstein:1999a}
    \textsc{Cannas~da Silva, A., Weinstein, A.:} \newblock
    \textit{Geometric Models for Noncommutative Algebras},
    Volume~10 in \textit{Berkeley Mathematics Lecture Notes}.
    AMS, 1999.

\bibitem {carinena.grabowski.marmo:2004a}
{\sc Cari{\~n}ena, J.~F., Grabowski, J., Marmo, G.: }\newblock {\em Courant
  algebroid and {L}ie bialgebroid contractions}.
\newblock J. Phys. A  {\bf 37}.19 (2004), 5189--5202.

\bibitem {courant:1990a}
{\sc Courant, T.~J.: }\newblock {\em {D}irac Manifolds}.
\newblock Trans. AMS  {\bf 319}.2 (1990), 631--661.

\bibitem {crainic.fernandes:2003a}
{\sc Crainic, M., Fernandes, R.~L.: }\newblock {\em Integrability of {L}ie
  brackets}.
\newblock Ann. of Math. (2)  {\bf 157}.2 (2003), 575--620.

\bibitem {crainic.moerdijk:2004a:pre}
{\sc Crainic, M., Moerdijk, I.: }\newblock {\em {D}eformations of {L}ie
  Brackets: Cohomologial Aspects}.
\newblock Preprint  {\bf math.DG/0403434} (2004).

\bibitem {dito.sternheimer:2002a}
{\sc Dito, G., Sternheimer, D.: }\newblock {\em Deformation quantization:
  genesis, developments and metamorphoses}.
\newblock In: {\sc Halbout, G. (eds.): }\newblock {\em Deformation
  quantization}, vol.~1 in {\em IRMA Lectures in Mathematics and Theoretical
  Physics},   9--54. Walter de Gruyter, Berlin, New York, 2002.

\bibitem {eilks:2004a}
{\sc Eilks, C.: }\newblock {\em {BRST}-{R}eduktion linearer {Z}wangsbedingungen
  im {R}ahmen der {D}eformationsquantisierung}.
\newblock master thesis, Fakult{\"{a}}t f{\"{u}}r Physik,
  Albert-Ludwigs-Universit{\"{a}}t, Freiburg, 2004.

\bibitem {fedosov:1996a}
{\sc Fedosov, B.~V.: }\newblock {\em Deformation Quantization and Index
  Theory}.
\newblock Akademie Verlag, Berlin, 1996.

\bibitem {gerstenhaber:1964a}
{\sc Gerstenhaber, M.: }\newblock {\em On the Deformation of Rings and
  Algebras}.
\newblock Ann. Math.  {\bf 79} (1964), 59--103.

\bibitem {gerstenhaber.schack:1988a}
{\sc Gerstenhaber, M., Schack, S.~D.: }\newblock {\em Algebraic Cohomology and
  Deformation Theory}.
\newblock In: {\sc Hazewinkel, M., Gerstenhaber, M. (eds.): }\newblock {\em
  Deformation Theory of Algebras and Structures and Applications},   13--264.
  Kluwer Academic Press, Dordrecht, 1988.

\bibitem {gualtieri:2003a}
{\sc Gualtieri, M.: }\newblock {\em Generalized complex geometry}.
\newblock PhD thesis, St John's College, University of Oxford, Oxford, 2003.
\newblock math.DG/0401221.

\bibitem {gutt:2000a}
{\sc Gutt, S.: }\newblock {\em Variations on deformation quantization}.
\newblock In: {\sc Dito, G., Sternheimer, D. (eds.): }\newblock {\em
  Conf{\'e}rence Mosh{\'e} Flato 1999. Quantization, Deformations, and
  Symmetries}, {\em Mathematical Physics Studies} no. {\bf 21},   217--254.
  Kluwer Academic Publishers, Dordrecht, Boston, London, 2000.

\bibitem {keller:2004a}
{\sc Keller, F.: }\newblock {\em Deformation von {L}ie-{A}lgebroiden und
  {D}irac-{S}trukturen}.
\newblock master thesis, Fakult{\"{a}}t f{\"{u}}r Mathematik und Physik,
  Physikalisches Institut, Albert-Ludwigs-Universit{\"{a}}t, Freiburg, 2004.

\bibitem {kontsevich:2003a}
{\sc Kontsevich, M.: }\newblock {\em Deformation Quantization of {P}oisson
  manifolds}.
\newblock Lett. Math. Phys.  {\bf 66} (2003), 157--216.

\bibitem {kosmann-schwarzbach:1996a}
{\sc Kosmann-Schwarzbach, Y.: }\newblock {\em From {P}oisson algebras to
  {G}erstenhaber algebras}.
\newblock Ann. Inst. Fourier (Grenoble)  {\bf 46}.5 (1996), 1243--1274.

\bibitem {kosmann-schwarzbach:2004b}
{\sc Kosmann-Schwarzbach, Y.: }\newblock {\em Derived brackets}.
\newblock Lett. Math. Phys.  {\bf 69} (2004), 61--87.

\bibitem {kosmann-schwarzbach:2005a}
{\sc Kosmann-Schwarzbach, Y.: }\newblock {\em Quasi, twisted, and all
  that{$\ldots$}in {P}oisson geometry and {L}ie algebroid theory}.
\newblock In: {\sc Marsden, J.~E., Ratiu, T.~S. (eds.): }\newblock {\em The
  breadth of symplectic and Poisson geometry}, vol. 232 in {\em Progress in
  Mathematics},   363--389. Birkh\"auser Boston Inc., Boston, MA, 2005.
\newblock Festschrift in honor of Alan Weinstein.

\bibitem {kosmann-schwarzbach.magri:1990a}
{\sc Kosmann-Schwarzbach, Y., Magri, F.: }\newblock {\em Poisson {N}ijenhuis
  structures}.
\newblock Ann. Inst. H. Poincar\'e Phys. Th\'eor.  {\bf 53} (1990), 35--81.

\bibitem {laurent-gengoux.ponte.xu:2005a}
{\sc Laurent-Gengoux, C., Ponte, D., Xu, P.: }\newblock {\em Universal lifting
  theorem and quasi-Poisson groupoids}.
\newblock math/0507396, 2005.

\bibitem {liu.weinstein.xu:1997a}
{\sc Liu, Z., Weinstein, A., Xu, P.: }\newblock {\em Manin triples for {L}ie
  bialgebroids}.
\newblock J. Differential Geom.  {\bf 45}.3 (1997), 547--574.

\bibitem {mackenzie:2005a}
{\sc Mackenzie, K. C.~H.: }\newblock {\em General Theory of {L}ie Groupoids and
  {L}ie Algebroids}, vol. 213 in {\em London Mathematical Society Lecture Note
  Series}.
\newblock Cambridge University Press, Cambridge, UK, 2005.

\bibitem {mackenzie.xu:1994a}
{\sc Mackenzie, K. C.~H., Xu, P.: }\newblock {\em Lie bialgebroids and
  {P}oisson groupoids}.
\newblock Duke Math. J.  {\bf 73}.2 (1994), 415--452.

\bibitem {marle:2002a}
{\sc Marle, C.-M.: }\newblock {\em Differential calculus on a {L}ie algebroid
  and {P}oisson manifolds}.
\newblock In: {\em The J. A. Pereira da Silva birthday schrift}, vol.~32 in
  {\em Textos Mat. S\'er. B},   83--149. Univ. Coimbra, Coimbra, 2002.

\bibitem {rothstein:1991a}
{\sc Rothstein, M.: }\newblock {\em The structure of supersymplectic
  supermanifolds}.
\newblock In: {\sc Bartocci, C., Bruzzo, U., Cianci, R. (eds.): }\newblock {\em
  Differential geometric methods in theoretical physics (Rapallo, 1990)}, vol.
  375 in {\em Lecture Notes in Physics},   331--343. Springer, Berlin, 1991.
\newblock Proceedings of the Nineteenth International Conference held in
  Rapallo, June 19--24, 1990.

\bibitem {roytenberg:1999a}
{\sc Roytenberg, D.: }\newblock {\em Courant Algebroids, derived brackets and
  even symplectic supermanifolds}.
\newblock PhD thesis, UC Berkeley, Berkeley, 1999.
\newblock math.DG/9910078.

\bibitem {roytenberg:2002b}
{\sc Roytenberg, D.: }\newblock {\em On the structure of graded symplectic
  supermanifolds and {C}ourant algebroids}.
\newblock In: {\sc Voronov, T. (eds.): }\newblock {\em Quantization, Poisson
  brackets and beyond (Manchester, 2001)}, vol. 315 in {\em Contemporary
  Mathematics},   169--185. American Mathematical Society, Providence, RI,
  2002.
\newblock Papers from the London Mathematical Society Regional Meeting held
  July 6, 2001 and the Workshop on Quantization, Deformations, and New
  Homological and Categorical Methods in Mathematical Physics held at the
  University of Manchester, Manchester, July 7--13, 2001.

\bibitem {roytenberg:2002a}
{\sc Roytenberg, D.: }\newblock {\em Quasi-{L}ie bialgebroids and twisted
  {P}oisson manifolds}.
\newblock Lett. Math. Phys.  {\bf 61}.2 (2002), 123--137.

\bibitem {severa:2005a:pre}
{\sc {\v{S}}evera, P.: }\newblock {\em On Deformation Quantization of Dirac
  Structures}.
\newblock Preprint  {\bf math.QA/0511403} (2005).

\bibitem {severa.weinstein:2001a}
{\sc {\v{S}}evera, P., Weinstein, A.: }\newblock {\em Poisson Geometry with a
  $3$-Form Background}.
\newblock In: {\sc Maeda, Y., Watamura, S. (eds.): }\newblock {\em
  Noncommutative Geometry and String Theory}, vol. 144 in {\em Prog. Theo.
  Phys. Suppl.},   145--154. Yukawa Institute for Theoretical Physics, 2001.
\newblock Proceedings of the International Workshop on Noncommutative Geometry
  and String Theory.

\bibitem {uchino:2002a}
{\sc Uchino, K.: }\newblock {\em Remarks on the definition of a {C}ourant
  algebroid}.
\newblock Lett. Math. Phys.  {\bf 60}.2 (2002), 171--175.

\bibitem {weinstein:1994a}
{\sc Weinstein, A.: }\newblock {\em Deformation Quantization}.
\newblock S\'eminaire Bourbaki 46\`eme ann\'ee  {\bf 789} (1994).

\end{thebibliography}

\end{footnotesize}

\end{document}